\documentclass[12pt]{article}
\usepackage[frenchb]{babel}
\usepackage{amssymb}
\usepackage{amsthm}

\newtheorem{defi}{D\'efinition}[section]
\newtheorem{Theoreme}[defi]{Th\'eor{\`e}me}
\newtheorem{lem}[defi]{Lemme}
\newtheorem{coro}[defi]{Corollaire}
\newtheorem{prop}[defi]{Proposition}

\theoremstyle{definition}
\newtheorem{rem}[defi]{Remarque}

\newcommand{\bK}{{\mathbb K}}
\newcommand{\bQ}{{\mathbb Q}}
\newcommand{\bR}{{\mathbb R}}
\newcommand{\bZ}{{\mathbb Z}}

\newcommand{\cA}{{\mathcal A}}
\newcommand{\cB}{{\mathcal B}}
\newcommand{\cD}{{\mathcal D}}
\newcommand{\cF}{{\mathcal F}}
\newcommand{\cG}{{\mathcal G}}
\newcommand{\cH}{{\mathcal H}}
\newcommand{\cM}{{\mathcal M}}
\newcommand{\cN}{{\mathcal N}}
\newcommand{\cP}{{\mathcal P}}

\newcommand{\hgt}{{\rm h}}
\newcommand{\length}{{\rm L}}

\title{Approximation diophantienne et approximants de Hermite-Pad\'e de type I
de fonctions exponentielles}
\author{Samy Kh\'emira et Paul Voutier}
\date{}

\begin{document}

\maketitle

\begin{abstract}
En utilisant des approximants de Hermite-Pad\'e de fonctions exponentielles,
ainsi que des d\'eterminants d'interpolation de Laurent, nous minorons la distance
entre un nombre alg\'ebrique et l'exponentielle d'un nombre alg\'ebrique non nul.
\end{abstract}

\begin{center}
{\bf Abstract}
\end{center}
\vspace{-3.0mm}
{\small We use Hermite-Pad\'e approximants of exponential functions along with
Laurent's interpolation determinants to obtain lower bounds for the distance
between an algebraic number and the exponential of another non-zero algebraic
number.}



\section{Introduction}

L'\'etude de la nature arithm\'etique de valeurs particuli\`eres de la fonction
exponentielle a longtemps \'et\'e l'un des th\`emes principaux de la th\'eorie
des nombres. Par exemple, en 1873, Hermite \cite{He1} a prouv\'e que $e$ est
transcendant. Peu apr\`es, Lindemann a d\'emontr\'e que $e^{\beta}$ est transcendant
si ${\beta}$ est un nombre alg\'ebrique diff\'erent de z\'ero. De cette fa\c{c}on,
il a \'et\'e capable de montrer que ${\pi}$ est transcendant, puisque $e^{i\pi}=-1$,
et ainsi de r\'esoudre l'ancien probl\`eme grec de la quadrature du cercle.

Quand on sait qu'un nombre est transcendant, il est pertinent de se poser la
question de son approximation par des nombres alg\'ebriques. De tels \'enonc\'es
quantitatifs sont connus sous le nom de ``mesures de transcendance''. Les plus
anciens r\'esultats dans cette direction remontent \`a 1899 quand Borel
\cite{Bo} a obtenu la premi\`ere mesure de transcendance de $e$. Nous renvoyons
le lecteur \`a \cite{FeNe} Chap.~2 pour plus d'informations sur l'histoire de
ce sujet, qui est le th\`eme de cet article.

\subsection{\'Enonc\'es}

Soit $\alpha$ un nombre alg\'ebrique de degr\'{e} $d$ sur $\bQ$, et donc le
polyn\^{o}me minimal sur $\bZ$ s'\'{e}crit $a\prod_{i=1}^{d} (X-\alpha^{(i)})$,
o\`{u} les racines $\alpha^{(i)}$ sont des nombres complexes, nous d\'esignons
par
$$
\hgt(\alpha) = \frac{1}{d} \left( \log |a| + \sum_{i=1}^{d} \log \max \left( 1, \left| \alpha^{(i)} \right| \right) \right),
$$
la hauteur logarithmique absolue de Weil du nombre alg\'ebrique $\alpha$.


Pour $m$ et $n$ des entiers naturels, nous posons
$$
D_{m,n}
= \frac{m!}{\displaystyle \prod_{
\stackrel{q \leq n}{q {\rm \; premier}}}q^{v_{q}(m!)}}
$$
et nous d\'efinissons $d_{n}$ comme le plus petit commun multiple de $1,\ldots,n$
avec $d_{0}=1$.

Nous allons d\'emontrer les r\'esultats \'enonc\'es ci-dessous.

\begin{Theoreme}
\label{theo:4.6bisalg}
Soient $\alpha$ et $\beta$ deux nombres alg\'ebriques avec $\beta \neq 0$.
Posons $D=[\bQ(\alpha,\beta):\bQ]/[\bR(\alpha,\beta):\bR]$. Soient $\cA \geq 1$
et $\cB \geq 1$ deux nombres r\'eels tels que
$D\hgt(\alpha)-\log(\max(1,|\alpha|)) \leq \log \cA$ et
$D\hgt(\beta)-\log(\max(1,|\beta|)) \leq \log \cB$.

Soient $K$ et $L$ des entiers strictement positifs avec $L \geq 2$ et $E>1$ un
nombre r\'eel tel que l'in\'egalit\'e suivante soit v\'erifi\'ee
\begin{eqnarray}
\label{eq:theoTalg}
KL\log E
& \geq & DKL\log 2 + D(K-1) \log \left( e\sqrt{3L}d_{L-1} \right)
         + D\log \left( D_{K-1,L-1} \right) \nonumber \\
&      & + D\log \left( (4e)^{L-1} \min \left( d_{L-2}^{K-1}, (L-2)! \right) \right)
         + \log ((K-1)!) \\
&      & + (K-1)\log (\cB/2) + (L-1)\log (\cA/2)
         + LE|\beta| + L\log E. \nonumber
\end{eqnarray}

Enfin soit $\epsilon = \left| \exp(\beta)-\alpha \right|$. Alors 
$$
\epsilon \geq E^{-KL}.
$$
\end{Theoreme}

\begin{rem}
L'in\'{e}galit\'{e} $L \geq 2$ n'est pas limit\'{e}. Le cas $L=1$ est sans
int\'{e}r\^{e}t, puisque les fonctions auxiliaires sont des polyn\^{o}mes.
\end{rem}

\begin{rem}
Sous les hypoth\`eses du th\'eor\`eme~\ref{theo:4.6bisalg}, le th\'eor\`eme
de Hermite-Lindemann \cite[Theorem~2.2, \S~2.3]{FeNe} et
\cite[Theorem~1.2, page 2]{Wa} s'\'enonce $\epsilon > 0$.
\end{rem} 

\begin{rem}
De (\ref{eq:theoTalg}), nous en d\'eduisons que $KL\log E \geq DKL\log 2$ et
$E \geq 2^{D}$. \`A partir de cette m\^{e}me in\'{e}galit\'e, nous obtenons
\'egale{\-}ment $KL\log E \geq LE|\beta|$; et en combinant cela avec l'in\'{e}galit\'{e}
pr\'{e}c\'{e}{\-}dente pour $E$, nous avons alors $KL\log E \geq L2^{D}|\beta|$.
Donc $KL\log E$ cro\^{i}t exponentiellement en $D$. Par la suite, le
th\'eor\`eme~\ref{theo:4.6bisalg} est le mieux adapt\'{e} pour les petites
valeurs de $D$.
\end{rem} 

Un cas particulier du th\'eor\`eme~\ref{theo:4.6bisalg} concerne
l'approximation par des entiers alg\'{e}briques de l'exponentielle d'un entier
alg\'{e}brique. Quand on se restreint aux corps imaginaires quadratiques, on en
d\'eduit:

\begin{coro}
\label{coro:4}
Quand $\alpha$ et $\beta$ sont deux entiers alg\'{e}briques dans le m\^{e}me
corps imaginaire quadratique avec $|\beta|$ suffisamment grand, on a
$$
\left| e^{\beta}-\alpha \right| \geq |\beta|^{-276.55 \ldots |\beta|}.
$$
\end{coro}

Historiquement, la premi\`ere estimation en direction de ce corollaire~\ref{coro:4}
est due \`a Mahler (voir \cite{Ma} et \cite{FeNe} page~105) avec une constante
\'egale \`a $33$ pour $\alpha, \beta \in \bZ$. Elle fut am\'elior\'ee et la
constante fut abaiss\'ee \`a $21$ (par Mignotte \cite{Mi}) et $19.187$ (meilleure
estimation connue \`a ce jour, par Wielonsky \cite{Wi}), encore dans le cas o\`u
$\alpha, \beta \in \bZ$.

De telles minorations ont des applications en informatique th\' eorique: dans
\cite{MuTi}, J.-M. Muller et A. Tisserand ont en effet utilis\'e les r\'esultats de
\cite{NeWa}.
Nous projetons, dans un futur article, de d\'{e}duire du th\'eor\`eme~\ref{theo:4.6bisalg}
des am\' elio{\-}rations des th\' eor\`emes 1--5 de \cite{NeWa} pour $D=1$.

Enfin, les auteurs souhaitent exprimer leur gratitude \`a ``les deux Michels''
(Michel Laurent et Michel Waldschmidt) pour leur aide pr\'{e}cieuse au cours de
ce travail. Elle a commenc\'{e} dans la recherche entreprise par le premier
auteur en travaillant avec Michel Laurent et qui est finalement devenu, pour sa
Th\`ese de Doctorat, l'\'el\`eve de Michel Waldschmidt. Les deux Michels ont
g\'{e}n\'{e}reusement offert de leur temps et leur expertise au cours de la
pr\'{e}paration de cet article et ont patiemment donn\'{e} des r\'{e}ponses
perspicaces aux nombreuses questions pos\'{e}es par le second auteur. Sans
leurs efforts, cet article n'aurait probablement jamais vu le jour sous cette
forme.

\subsection{Introduction d'un d\'eterminant et plan de la d\'emonstration du th\'eor\`eme~\ref{theo:4.6bisalg}}

Voici le sch\'ema de la d\'emonstration du th\'eor\`eme~\ref{theo:4.6bisalg}.\\

Pour $k$ et $\ell$ deux entiers naturels, consid\'erons les fonctions
$$
\Phi_{k,\ell}(z)=z^{k}e^{\ell z }.
$$

Pour $K$ et $L$ des entiers strictement positifs, posons $S=KL$ et consid\'erons la
matrice $\cM(x,y) \in {\mathbb M}_{S,S}(\bQ[x,y])$ d\'efinie par
$$
\cM(x,y)=
\left(
	\begin{array}{l}
		\cM_{0} \\
		\cM_{1}(x,y) 
	\end{array}
\right)
$$
avec
$$
\cM_{0}
= \left( \left( \frac{{\rm d}}{{\rm d}z} \right)^{s} \left( \Phi_{k,\ell}(z) \right)\left(0\right)
  \right)_{\stackrel{0 \leq s \leq S-2}{0 \leq k \leq K-1, 0 \leq \ell \leq L-1}},
$$
o\`{u} $s$ d\'{e}signe l'indice de ligne et o\`{u} les couples $(k,\ell)$
param\'{e}trisent les colonnes, et
\begin{eqnarray*}
\cM_{1}(x,y)
& = & \left( \delta^{\mu}\left(X^{k}Y^{\ell}\right)(x,y)
      \right)_{0 \leq k \leq K-1, 0 \leq \ell \leq L-1} \\
& = & \left( y^{\ell}\sum_{j=0}^{\min(\mu, k)}
         {\mu \choose j} \frac{k!}{(k-j)!}x^{k-j}\ell^{\mu-j}
      \right)_{0 \leq k \leq K-1, 0 \leq \ell \leq L-1},
\end{eqnarray*}
o\`u $\mu$ d\'{e}signe un entier naturel qui sera choisi plus tard, $0^0=1$,
et $\delta$ est l'op\'erateur de d\'erivation d\'efini par:
\begin{equation}
\label{eq:op}
\delta= \frac{\partial}{\partial X}+Y\frac{\partial}{\partial Y}.
\end{equation}

\begin{rem}
\label{rem:4.13}
On a
$$
\left( \frac{{\rm d}}{{\rm d}z} \right)^{s} \left( \Phi_{k,\ell}(z) \right)\left(0\right)
= \left\{ \begin{array}{ll}
                    0  & \mbox{si $k > s$} \\
                    \displaystyle {s \choose k} k! \ell^{s-k} & \mbox{sinon.}
                    \end{array}
            \right.
$$
\end{rem}

\begin{defi}
\label{defi:bigF}
Nous posons tout d'abord $F(x,y) = \det \cM(x,y)$, et ensuite
$\cD = F(\beta,\alpha) = \det \left( \cM \left( \beta, \alpha \right) \right)$.
\end{defi}

La remarque~\ref{rem:4.13} et la d\'efinition de $\cM_{1}(x,y)$ permettent d'en
d\'eduire que $F(x,y)\in \bZ[x,y]$.

\begin{defi}
\label{defi:wkl}
Pour\ $0 \leq k \leq K-1$, $0 \leq \ell \leq L-1$, d\'efinissons les nom{\-}bres
complexes $w_{k,\ell}$ par:
$$
w_{k,\ell} = \frac{\alpha^{\ell}-e^{\beta\ell}}{|e^{\beta}-\alpha|}
\sum_{j=0}^{\min (\mu, k)} {\mu \choose j}
\frac{k!\ell^{\mu-j}}{(k-j)!} \beta^{k-j}.
$$
\end{defi}

\begin{rem}
Notons que l'on a 
$$
\delta^{\mu}\left(X^{k}Y^{\ell}\right)(\beta, \alpha)
=\Phi_{k,\ell}^{(\mu)}(\beta) + \epsilon w_{k,\ell}
$$
et
$$
\cD =
\det
\left(
	\begin{array}{l}
		\cM_{0} \\
		\left( \Phi_{k,\ell}^{(\mu)}(\beta)+\epsilon w_{k,\ell}
		\right)_{0 \leq k \leq K-1,0 \leq \ell \leq L-1}
	\end{array}
\right).
$$
\end{rem}

\begin{rem}
D'autres choix sont possibles pour $\cM(x,y)$. Par ex{\-}emple, on pourrait
utiliser une matrice $\cM_{1}(x,y)$ comportant $T$ lignes index\'{e}es par
plusieurs d\'{e}rivations $\mu_{1}$, \ldots, $\mu_{T}$, ou m\^{e}me consid\'{e}rer
plusieurs points multiples entiers de $\beta$. Mais les travaux du premier
auteur, dans la Section 3.5 de sa th\`{e}se \cite{Kh}, ainsi que des travaux
informatiques effectu\'es par le second auteur, sugg\`{e}rent (mais ne
prouvent pas!) que le choix $T=1$ d\'evelopp\'e ici pr\'esent\'e pourrait
\^{e}tre le meilleur.
\end{rem}

La m\'ethode habituelle en th\'eorie des nombres transcendants
(quand on utilise les d\'eterminants d'interpolation de Laurent) consiste \`a
montrer que $\cD$ n'est pas nul (lemme de z\'eros), \`a le majorer par
un argument analytique (lemme de Schwarz), \`a le minorer par un argument
arithm\'etique (in\'egalit\'e de Liouville) et \`a comparer les deux estimations
pour obtenir la conclusion. Nous suivrons cette d\'emarche, mais le nombre qui
interviendra n'est pas $\cD$ lui-m\^eme, c'est un nombre qui lui est \'etroitement
li\'e (il s'agit de $\cG_{\beta,\alpha}$ donn\'e par la formule (\ref{eq:defG})
plus loin). Le lemme de z\'eros est la proposition~\ref{prop:NW}, la minoration
arithm\'etique est la proposition~\ref{prop:trivalg} et la majoration analytique
est la proposition~\ref{prop:4.17alg} ci-dessous.

Un point important de cette \'etude se trouve dans le lemme de z\'eros \'enonc\'e
dans cet article, qui est un raffinement de celui utilis\'e dans \cite{NeWa}.

\section{Lemmes auxiliaires}

\subsection{Un nouveau lemme de z\'eros}

Nous allons \'enoncer et d\'emontrer un lemme de z\'eros qui pr\'ecise celui de
Yu. V. Nesterenko et M. Waldschmidt dans \cite{NeWa}.

Soient $\bK$ un corps de caract\'eristique nulle et $\delta$ l'op\'erateur d\'efini
en (\ref{eq:op}) sur l'espace $\bK[X,Y]$ des polyn\^omes en deux variables.

\begin{prop}
\label{prop:NW}
Soient $M$ un entier positif et $D_{0}$, $D_{1}$, $S_{1}, \ldots$, $S_{M}$
des entiers naturels v\'erifiant
\begin{equation}
\label{eq:hyp}
S_{1} + \cdots + S_{M} > (D_{0}+M)(D_{1}+1)-M.
\end{equation}
Soient $(\zeta_{1}, \eta_{1}), \ldots, (\zeta_{M}, \eta_{M})$ des couples de
$\bK \times \bK^{*}$ o\`u $\zeta_{1}, \ldots, \zeta_{M}$ sont
deux \`a deux distincts. Il n'existe alors pas de polyn\^ome non nul
$P \in \bK[X,Y]$, de degr\'e $\le D_{0}$ en $X$ et de degr\'e
$\le D_{1}$ en $Y$, qui satisfasse
\begin{equation}
\label{eq:cond}
\delta^{\sigma} P(\zeta_{\kappa}, \eta_{\kappa}) = 0
\end{equation}
pour tout couple $\sigma, \kappa$ v\'{e}rifiant
$0 \leq \sigma <S_{\kappa}$ et $1 \leq \kappa \leq M$.
\end{prop}

\begin{dem}
Supposons qu'il existe un polyn\^ome non nul $P$ de degr\'e $\leq D_{0}$ en $X$
et de degr\'e $\leq D_{1}$ en $Y$, qui satisfasse les \'egalit\'es (\ref{eq:cond}).
Nous allons montrer que la condition (\ref{eq:hyp}) n'est pas satisfaite.

Nous pouvons supposer, sans perte de g\'en\'eralit\'e, que $Y$ ne divise pas
le polyn\^ome $P$ (puisque $\eta_{\kappa} \neq 0$), et que $P$ a un degr\'e
sup\'erieur ou \'egal \`a $1$ par rapport \`a la variable $Y$ (le cas o\`u $P$
est de degr\'e $0$ en la variable $Y$ est celui, trivial, d'un polyn\^ome en une
variable). D\'efinissons les entiers $n$, $k_{0}, \ldots, k_{n}$,
$m_{0}, \ldots, m_{n}$, les polyn\^omes $Q_{i} \in \bK[X]$
($0 \leq i \leq n$) et les \'el\'ements $b_{0},\ldots,b_{n}$ de $\bK^{*}$
par les conditions suivantes:
\begin{eqnarray*}
k_{0}    & = & 0 < k_{1} < \cdots < k_{n} \leq D_{1}, \\
P(X,Y)   & = & \sum_{i=0}^{n} Q_{i}(X)Y^{k_{i}}, \\
Q_{i}(X) & = & b_{i}X^{m_{i}} +\cdots \in \bK[X],
\hspace{3.0mm} b_{i} \in \bK^{*} \hspace{3.0mm} (i=0,\ldots,n).
\end{eqnarray*}

Pour $0 \leq \sigma \leq n$, nous d\'efinissons les polyn\^omes
$Q_{\sigma i} \in \bK[X]$ par
\begin{equation}
\label{eq:cond2}
\delta^{\sigma}P(X,Y)=\sum_{i=0}^{n}Q_{\sigma i}(X)Y^{k_{i}},
\end{equation}
de sorte que
$$
Q_{\sigma i}(X)
= \sum_{j=0}^{\sigma} {\sigma \choose j} Q_{i}^{(\sigma-j)}(X) k_{i}^{j}
= b_{i} k_{i}^{\sigma}X^{m_{i}} + \cdots.
$$
On pose
\begin{eqnarray*}
\Delta(X) &= & \det \left( Q_{\sigma i}(X) \right)_{0 \leq i,\sigma \leq n}
=\det \left( b_{i}k_{i}^{\sigma}X^{m_{i}} + \cdots \right)_{0 \leq i,\sigma \leq n} \\
& = & b_{0} \cdots b_{n} B X^{m_{0}+\cdots+m_{n}}+\cdots,
\end{eqnarray*}
o\`u $B$ est le d\'eterminant de Vandermonde construit \`a partir des nombres
$k_{0}, \ldots, k_{n}$, donc $B \neq 0$. De (\ref{eq:cond2}),
nous d\'eduisons, par les formules de Cramer,
$$
\Delta(X) = \sum_{\sigma=0}^{n} \Delta_{\sigma}(X,Y) \delta^{\sigma}P(X,Y)
\hspace{3.0mm} {\rm avec} \hspace{3.0mm}
\Delta_{\sigma}(X,Y) \in \bK[X,Y].
$$
Par cons\'equent, pour $1 \leq j \leq M$ et $0 \leq \tau < \max \left\{ 0,S _{j}-n \right\}$,
on peut \'ecrire
$$
\Delta^{(\tau)}(\zeta_{j})
= \sum_{\sigma=0}^{n+\tau} c_{\tau,j,\sigma} \delta^{\sigma}P(\zeta_{j},\eta_{j})
= 0
$$
avec $c_{\tau_{j},j,\sigma} \in \bK$. \'Etant donn\'e que les nombres r\'eels
$\zeta_{j}$ sont deux \`a deux distincts, par hypoth\`ese, et puisque la somme
des multiplicit\'es des z\'eros d'un polyn\^ome en une variable est \'egale au
degr\'e de ce polyn\^ome, on obtient
$$
\sum_{j=1}^{M}(S _{j}-n)
\leq \sum_{j=1}^{M} \max \left\{ 0,S _{j}-n \right\}
\leq \deg \Delta(X).
$$

Or $n \leq D_{1}$ et
$\deg \Delta(X) = m_{0} + \cdots + m_{n} \leq (n+1)D_{0} \leq (D_{1}+1)D_{0}$. Donc
$$
\sum_{j=1}^{M} S_{j} \leq (D_{1}+1)D_{0} + nM \leq (D_{1}+1)D_{0}+D_{1}M
= (D_{0}+M)(D_{1}+1)-M.
$$

Ceci montre que la condition (\ref{eq:hyp}) n'est pas satisfaite.
La proposition~\ref{prop:NW} en r\'esulte. \hfill $\Box$
\end{dem}

\begin{rem}
Le lemme~$2$ de \cite{NeWa} est le cas particulier de la proposition~$\ref{prop:NW}$
dans lequel $S_{1} = \cdots = S_{M}$.
\end{rem}

\begin{rem}
Un exemple avec $D_{0}=0$ o\`u le lemme de z\'eros pr\'ec\'edent
$($proposition~\ref{prop:NW}$)$ est optimal est donn\'e par le polyn\^ome
$\cP(X,Y)=(Y-1)^{D_{1}}$ qui v\'erifie $(\ref{eq:cond})$ avec
$$
S_{1} = \cdots = S_{M} = D_{1},
$$
pour les points $(\zeta_{\mu},\eta_{\mu})=(\mu,1)$, $1 \leq \mu \leq M$.
\end{rem}

\subsubsection{Une application}

Consid\'erons le polyn\^ome $\cH(X,Y) \in \bZ[X,Y]$ d\'efini par
\begin{equation}
\label{eq:bigH}
\cH(X,Y) = \det \left(
\begin{array}{l}
    \cM_{0} \\
    \left( X^{k}Y^{\ell} \right)_{0 \leq k \leq K-1,0 \leq \ell \leq L-1}
\end{array} \right).
\end{equation}

\begin{lem}
\label{lem:h-not-zero}
Soient $\alpha$ et $\beta$ deux nombres alg\'ebriques avec $\alpha, \beta \neq 0$
et $L \geq 2$ un entier. Alors il existe un entier naturel $\mu \leq L-2$ tel
que $\delta^{\mu} \cH(\beta, \alpha) \neq 0$.
\end{lem}

\begin{dem}
Le polyn\^ome $\cH(X,Y)$ est de degr\'e au plus $K-1$ en $X$ et au plus $L-1$
en $Y$.

De plus,
\begin{eqnarray*}
\delta^{\mu} \cH(z, e^{z})
& = & \det \left(
\begin{array}{l}
    \cM_{0} \\
    \left( \delta^{\mu} \left( X^{k}Y^{\ell} \right) (z, e^{z}) \right)_{0 \leq k \leq K-1,0 \leq \ell \leq L-1}
\end{array} \right) \\
& = & \det \left(
\begin{array}{l}
    \cM_{0} \\
    \left( \Phi_{k,\ell}^{(\mu)} (z) \right)_{0 \leq k \leq K-1,0 \leq \ell \leq L-1}
\end{array} \right),
\end{eqnarray*}
d'apr\`es la d\'efinition des $\Phi_{k,\ell}(z)$.

Notons que $S=KL \geq 2$, car $L \geq 2$.

Dans le cas particulier o\`u $z=0$, pour tout entier naturel $\mu$ v\'erifiant
$\mu \leq S-2$, la derni\`ere ligne de la matrice est la m\^eme que la ligne,
dans $\cM_{0}$, index\'ee par $s=\mu$. Par cons\'equent $\delta^{\mu} \cH(0,1)=0$
pour tout entier naturel $\mu$ v\'erifiant $\mu <S-1$.

On applique la proposition~\ref{prop:NW} avec $D_{0}=K-1$, $D_{1}=L-1$, $M=2$,
$S_{1}=S-1=KL-1$, $(\zeta_{1}, \eta_{1})=(0,1)$ et $(\zeta_{2}, \eta_{2})=(\beta,\alpha)$.

Il est \'{e}galement n\'{e}cessaire, pour utiliser la proposition~\ref{prop:NW},
que $\cH(X,Y)$ ne soit pas le polyn\^{o}me nul. Dans la preuve de la
proposition~\ref{prop:3.10} nous avons montr\'{e} que la matrice $\cM_{0}$ est
de rang $S-1$. Nous allons en d\'eduire qu'au moins une des sous-matrices obtenue
en supprimant une des colonnes est de d\'eterminant non nul. Pour cela, supposons
que l'indice d'une telle sous-matrice est $(k_{0}, \ell_{0})$. La derni\`{e}re
ligne de la matrice dont le d\'{e}terminant est $\cH(X,Y)$ contient des
mon\^{o}mes distincts de la forme $X^{k}Y^{\ell}$; quand on d\'eveloppe
$\cH(X,Y)$ comme un d\'eterminant par rapport \`a la  derni\`ere ligne, le
polyn\^ome que l'on obtient est
$c_{k_{0}, \ell_{0}}X^{k_{0}}Y^{\ell_{0}}$ plus d'autres
termes de degr\'e strictement inf\'erieur et o\`{u} $c_{k_{0}, \ell_{0}} \neq 0$,
et donc $\cH(X,Y)$ n'est pas le polyn\^ome nul.

Ainsi, toutes les conditions sont r\'eunies pour appliquer notre lemme de z\'eros
(proposition~\ref{prop:NW}): on a
$$
S_{1} + S_{2} = KL-1+S_{2} \leq (D_{0}+2)(D_{1}+1)-2 = (K+1)L-2
$$
o\`u $S_{2} \leq L-1$.

Par cons\'equent il existe un entier $0 \leq \mu \leq L-2$ tel que
$\delta^{\mu}\cH(\beta,\alpha) \neq 0$ et le lemme~\ref{lem:h-not-zero}
r\'esulte de la d\'efinition~\ref{defi:bigF}
et de l'\'egalit\'e (\ref{eq:defG}). \hfill$\Box$
\end{dem}

\begin{defi}
\label{defi:mu}
On prend pour $\mu$ le plus petit tel entier.
\end{defi}

\`A partir de maintenant, $\mu$ d\'esignera cette valeur.

\subsection{Polyn\^{o}mes de Fel'dman}

Nous introduisons des polyn\^omes de Fel'dman afin d'affiner de futures estimations.

\begin{defi}
\label{defi:polfel}
Pour $\nu$ entier naturel, on d\'efinit le polyn\^ome de Fel'dman
$\cF_{\nu}$ d'indice $\nu$
$($et ses coefficients $\lambda_{j,\nu} \in \bQ$, $j=0, \ldots, \nu)$
par
$$
\cF_{\nu}(z) = \frac{z(z-1) \cdots (z-\nu+1)}{\nu!}
= \sum_{j=0}^{\nu}\lambda_{j,\nu}z^{j}.
$$
\end{defi}

Notons que $\lambda_{j,\nu}=s(\nu,j)/\nu!$, o\`u $s(\nu, j)$ est le $(\nu, j)$-\`eme
nombre de Stirling de premi\`ere esp\`ece, voir \cite[Section 24.1.3, p.824]{AbSt}.

\begin{lem}
\label{lem:fel2}
Pour tout $\nu \geq 0$, on a
$$
\sum_{j=0}^{\nu} j! \left| \lambda_{j,\nu} \right|
\leq 2^{\nu}.
$$
\end{lem}

\begin{dem}
Comme $\displaystyle \sum_{j=0}^{0} j! \left| \lambda_{j,0} \right|=1 = 2^{0}$
et $\displaystyle \sum_{j=0}^{1} j! \left| \lambda_{j,1} \right| =1 < 2^{1}$,
on peut supposer $\nu \geq 2$ et proc\'eder par r\'ecurrence sur $\nu$.

\'Etant donn\'e que $\lambda_{0,\nu}=0$ et $\nu!\lambda_{\nu,\nu}=1$ pour
$\nu \geq 2$, on peut \'ecrire
$$
\sum_{j=0}^{\nu} j! \left| \lambda_{j,\nu} \right|
= 1 + \sum_{j=1}^{\nu-1} j! \left| \lambda_{j,\nu} \right|.
$$

D'apr\`es \cite[Recurrences, Section 24.1.3.II.A, p.824]{AbSt} (aussi de la
relation $\cF_{\nu}(z)=((z-\nu+1)/\nu)\cF_{\nu-1}(z)$), on a 
\begin{eqnarray*}
j! \left|\lambda_{j,\nu}\right|
& = & \frac{j}{\nu}(j-1)!\left|\lambda_{j-1,\nu-1}\right|
+ \frac{\nu-1}{\nu}j!\left|\lambda_{j,\nu-1}\right|.
\end{eqnarray*}

Par cons\'equent
\begin{eqnarray*}
\sum_{j=0}^{\nu} j! \left| \lambda_{j,\nu} \right|
& = & 1 + \sum_{j=1}^{\nu-1} \frac{j}{\nu} (j-1)!\left|\lambda_{j-1,\nu-1}\right|
+ \frac{\nu-1}{\nu} \sum_{j=1}^{\nu-1} j!\left|\lambda_{j,\nu-1}\right| \\
& = & \sum_{j=0}^{\nu-1} \frac{j+1}{\nu} j!\left|\lambda_{j,\nu-1}\right|
+ \frac{\nu-1}{\nu} \sum_{j=1}^{\nu-1} j!\left|\lambda_{j,\nu-1}\right| \\
& \leq & \frac{2 \nu-1}{\nu} \sum_{j=0}^{\nu-1} j!\left|\lambda_{j,\nu-1}\right|
<2 \sum_{j=0}^{\nu-1} j!\left|\lambda_{j,\nu-1}\right|.
\end{eqnarray*}

L'hypoth\`ese de r\'ecurrence permet de conclure la d\'emonstration du
lem{\-}me~\ref{lem:fel2}.
\hfill$\Box$
\end{dem}

\begin{lem}
\label{lem:fel3}
Soient $k$ et $\nu$ deux entiers naturels. Pour tout entier $\ell$
et tout entier $u$ dans l'intervalle $0 \leq u \leq k$, on a
$$
d_{\nu}^{k}\cF^{(u)}_{\nu}(\ell) \in \bZ.
$$
\end{lem}

\begin{dem}
Nous renvoyons au Lemme 4 (paragraphe $4$) de \cite{NeWa} pour une
d\'e{\-}monstration de ce lemme (en observant que le polyn\^{o}me $\cF_{\nu}(z)$
est not\'{e} $\Delta(z,\nu,\nu)$ dans \cite{NeWa}).
\hfill$\Box$
\end{dem}

\subsection{Approximants de Hermite-Pad\'e de type I de fonctions exponentielles}

Nous allons d'abord d\'efinir la notion d'approximants de Hermite-Pad\'e de
fonctions analytiques, puis donner quelques propri\'et\'es des approximants
de Hermite-Pad\'e de fonctions exponentielles.

\subsubsection{D\'efinitions et propri\'et\'es}

\begin{defi}
\label{defi:1.1}
Soit $m$ un entier naturel non nul.

Soient $n_{0}, \ldots, n_{m}$ des entiers positifs, $f_{0},\ldots, f_{m}$ des
fonctions analytiques en $0$; on appelle {\it syst\`eme d'approximants de Pad\'e}
$($ou {\it approximants de Hermite-Pad\'e}$)$ {\it de type I} pour les fonctions
$f_{0}, \ldots, f_{m}$ et les param\`etres $n_{0},\ldots, n_{m}$, tout
$(m+1)$-uplet $(P_{0}, \ldots, P_{m})$ de polyn\^omes non tous nuls tel que
$$
\deg P_{\ell}(z) \leq n_{\ell}-1
\hspace{3.0mm} \mbox{et} \hspace{3.0mm}
{\rm ord}_{z=0} \left( \sum_{\ell=0}^{m}P_{\ell}(z)f_{\ell}(z) \right)
\geq \sigma - 1,
$$
o\`u $\sigma = \sum_{\ell=0}^{m}n_{\ell}$.
\end{defi}

Hermite (voir \cite{He1}, \cite{He2} et \cite{Ma}) a donn\'e les formules
suivantes pour les approximants de Hermite-Pad\'e de type I pour les fonctions
$e^{x_{0}z}, \ldots, e^{x_{m}z}$ et les param\`etres $n_{0}, \ldots, n_{m}$:

\begin{lem}
\label{lem:2.1}
Soient $C_{0}$ et $C_{\infty}$ deux cercles, le premier centr\'e en $0$
et de rayon strictement inf\'erieur \`a la plus petite distance s\'eparant 
deux $x_{j}$, et le second centr\'e en $0$ contenant tous les points
$x_{0}, \ldots, x_{m}$.

Soient
$$
P_{\ell}(z) = \frac{1}{2\pi i} \int_{C_{0}}
\frac{e^{z\zeta}}{\displaystyle \prod_{j=0}^{m}
\left(\zeta+x_{\ell}-x_{j}\right)^{n_{j}}} \rm{d}\zeta
$$
et
$$
R(z) = \frac{1}{2\pi i} \int_{C_{\infty}}
\frac{e^{z\zeta}}{\displaystyle \prod_{j=0}^{m}
\left(\zeta-x_{j}\right)^{n_{j}}} \rm{d}\zeta.
$$

Alors\\
$(1)$ $P_{0}, \ldots, P_{m}$ sont des polyn\^omes de degr\'e $n_{\ell}-1$,\\
$(2)$ $\displaystyle R(z)= \sum_{\ell=0}^{m}P_{\ell}(z)e^{x_{\ell}z}$ et\\
$(3)$ $\displaystyle {\rm ord}_{z=0} \left( R(z) \right) \geq \sigma - 1$
o\`u $\displaystyle \sigma = \sum_{\ell=0}^{m}n_{\ell}$.

Les polyn\^omes $P_{0}, \ldots, P_{m}$ sont donc des approximants de Hermite-Pad\'e
de type I pour les fonctions $e^{x_{0}z}, \ldots, e^{x_{m}z}$ et les param\`etres
$n_{0}, \ldots, n_{m}$.
\end{lem}

\begin{defi}
\label{def:Gamma}
Pour $0 \leq \ell \leq m$ et $0 \leq k \leq n_{\ell}-1$, nous d\'efinissons
$\Lambda(\ell,k)$ $($qui d\'epend de $n_{0}, \ldots, n_{m})$ comme l'ensemble
des multi-indices entiers
$\gamma=\left( \gamma_{p}\right)_{\stackrel{0 \leq p \leq m}{p \neq \ell}}
\in \bZ^{m}$,
avec $\gamma_{p} \geq 0$ pour $p$ compris entre $0$ et $m$, $p \neq \ell$, tels
que $\displaystyle\sum_{\stackrel{p=0}{p \neq \ell}}^{m}\gamma_{p}=n_{\ell}-k-1$.
\end{defi}

Sous les hypoth\`eses du lemme $\ref{lem:2.1}$, comme chaque $P_{\ell}(z)$ est
un polyn\^ome de degr\'e $n_{\ell}-1$, notons
$$
p_{\ell,k} = P_{\ell}^{(k)}(0), \hspace{3.0mm}
\mbox{$(0 \leq k \leq n_{\ell}-1, 0 \leq \ell \leq m)$}
$$
les d\'{e}riv\'{e}es \`{a} l'origine des polyn\^{o}mes $P_{\ell}$. On a donc
le d\'{e}veloppment
$$
P_{\ell}(z) = \sum_{k=0}^{n_{\ell}-1} p_{\ell,k} \frac{z^{k}}{k!}.
$$

\begin{lem}
\label{lem:2.2}
On a 
$$
p_{\ell,k} = \sum_{\gamma \in \Lambda(\ell,k)}
\prod_{\stackrel{p=0}{p \neq \ell}}^{m}
\left( \frac{(-1)^{\gamma_{p}}}{(x_{\ell}-x_{p})^{\gamma_{p}+n_{p}}}
{\gamma_{p}+n_{p}-1 \choose n_{p}-1} \right).
$$
\end{lem}

En particulier, le coefficient dominant
$\displaystyle \frac{1}{(n_{\ell}-1)!} p_{\ell,n_{\ell}-1}$
du polyn\^ome $P_{\ell}(z)$ est
$$
\frac{1}{(n_{\ell}-1)!}
\prod_{\stackrel{p=0}{p \neq \ell}}^{m}
\frac{1}{\left( x_{\ell}-x_{p} \right)^{n_{p}}}.
$$

\begin{dem}
Nous allons utiliser les formules int\'egrales du lemme~\ref{lem:2.1}.

En d\'eveloppant $ e^{z\zeta} $ en s\'erie enti{\`e}re, on a 
$\displaystyle
e^{z\zeta}=\sum_{j=0}^{+\infty} \frac{\left(z\zeta\right)^{j}}{j!}$
et cette s\'erie converge normalement pour $\zeta$ d\'ecrivant $C_0$. Du
lemme~\ref{lem:2.1} on d\'eduit
$$
P_{\ell}(z) = \sum_{k=0}^{n_{\ell}-1} \left( \frac{1}{2\pi i}
\int_{C_{0}} f_{\ell,k}(\zeta) {\rm d}\zeta\right) \frac{z^{k}}{k!}
$$
o\`u on a not\'e
$$
f_{\ell,k}(\zeta)
= \frac{\zeta^{k-n_{\ell}}}{\displaystyle \prod_{\stackrel{p=0}{p \neq \ell}}^{m}
\left(\zeta+x_{\ell}-x_{p}\right)^{n_{p}}}.
$$

En utilisant l'expression pr\'ec\'edente des $P_{\ell}(z)$, on obtient
\begin{equation}
p_{\ell,k} = \frac{1}{2\pi i} \int_{C_{0}}f_{\ell,k}(\zeta)\rm{d}\zeta.
\end{equation}

Les p{\^o}les \'eventuels de la fonction $f_{\ell,k}$ sont en $0$ et en $x_{p}-x_{\ell}$
($p=0, \ldots, m$ et $p \neq \ell$); seul le point z\'ero se trouve {\`a}
l'int\'erieur du cercle $C_{0}$. Nous allons \'etudier le d\'eveloppement
de cette fonction $f_{\ell,k}$ au voisinage de z\'ero, afin d'en d\'eterminer
le r\'esidu.

En d\'erivant $j$ fois la s\'erie g\'eom\'etrique, on obtient la formule
$$
\frac{j!}{(1+x)^{j+1}} = \sum_{n \geq j} n(n-1) \cdots (n-j+1)(-x)^{n-j}
$$
qui fournit l'\'egalit\'e
$$
\frac{1}{\left(\zeta+x_{\ell}-x_{p}\right)^{n_{p}}}
= \frac{1}{(n_{p}-1)!}
\sum_{n \geq n_{p}-1} \frac{(-1)^{n-n_{p}+1}}{(x_{\ell}-x_{p})^{n+1}}
\frac{n!}{(n-n_{p}+1)!} \zeta^{n-n_{p}+1}.
$$

Par un simple changement d'indices, on trouve
$$
\frac{1}{\left(\zeta+x_{\ell}-x_{p}\right)^{n_{p}}}
= \frac{1}{(n_{p}-1)!}
\sum_{n \geq 0} \frac{(-1)^{n}}{(x_{\ell}-x_{p})^{n+n_{p}}}
\frac{(n+n_{p}-1)!}{n!} \zeta^{n}.
$$

Le produit de telles expressions donne:
\begin{eqnarray*}
\frac{1}{\displaystyle \prod_{\stackrel{p=0}{p \neq \ell}}^{m}
\left( \zeta+x_{\ell}-x_{p} \right)^{n_{p}}}
& = & \prod_{\stackrel{p=0}{p \neq i}}^{m}
\left( \frac{1}{(n_{p}-1)!}
\sum_{n \geq 0} \frac{(-1)^{n}}{(x_{\ell}-x_{p})^{n+n_{p}}}
\frac{(n+n_{p}-1)!}{n!} \zeta^{n}\right) \\
& = & \left( \prod_{\stackrel{p=0}{p \neq \ell}}^{m} \frac{1}{(n_{p}-1)!} \right)
\prod_{\stackrel{p=0}{p \neq \ell}}^{m}
\left( \sum_{n \geq 0} \frac{(-1)^{n}}{n!}
\frac{(n+n_{p}-1)!}{(x_{\ell}-x_{p})^{n+n_{p}}} \zeta^{n}\right).
\end{eqnarray*}

Le coefficient du mon\^ome $\zeta^{n_{\ell}-k-1}$ dans le d\'eveloppement de
$$
\prod_{\stackrel{p=0}{p \neq \ell}}^{m}
\left( \zeta+x_{\ell}-x_{p} \right)^{-n_{p}},
$$
qui permettra d'obtenir le r\'esidu cherch\'e, vaut:
\begin{eqnarray*}
& & \left( \prod_{\stackrel{p=0}{p \neq \ell}}^{m}
\frac{1}{(n_{p}-1)!} \right)
\sum_{\Lambda(\ell,k)}
\prod_{\stackrel{p=0}{p \neq \ell}}^{m}
\left({(-1)^{\gamma_{p}}\over{(x_{\ell}-x_{p})^{\gamma_{p}+n_{p}}}}
{(\gamma_{p}+n_{p}-1)! \over \gamma_{p}!}\right) \\
& = & \sum_{\Lambda(\ell,k)}
\prod_{\stackrel{p=0}{p \neq \ell}}^{m}
\left({(-1)^{\gamma_{p}}\over{(x_{\ell}-x_{p})^{\gamma_{p}+n_{p}}}}
{\gamma_{p}+n_{p}-1 \choose n_{p}-1} \right).
\end{eqnarray*}

D'apr{\`e}s le th\'eor{\`e}me des r\'esidus, on en d\'eduit 
$$
p_{\ell,k}=\sum_{\Lambda(\ell,k)}
\prod_{\stackrel{p=0}{p \neq \ell}}^{m}
\left( {(-1)^{\gamma_{p}} \over (x_{\ell}-x_{p})^{\gamma_{p}+n_{p}}}
{\gamma_{p}+n_{p}-1 \choose n_{p}-1} \right).
$$

Le lemme~\ref{lem:2.2} est ainsi d\'emontr\'e.
\hfill$\Box$
\end{dem}

\subsubsection{D\'eterminant de Vandermonde g\'en\'eralis\'e}

\begin{defi}
\label{defi:2.4}
Soit $\bK$ un corps de caract\'eristique nulle. Quand $m$ est un entier naturel,
$n_{0}$, $\ldots$, $n_{m}$ des entiers strictement positifs et
$x_{0}, \ldots, x_{m}$ des \'el\'ements de $\bK$, on appelle
{\og d\'eterminant de Vandermonde g\'en\'eralis\'e\fg} le d\'eterminant
$$
\Delta = \left| B_{0} \ldots B_{m} \right|,
$$
o{\`u}, pour $0\le \ell\le m$, $B_{\ell}$ d\'esigne la matrice {\`a}
$\displaystyle \sigma=\sum_{i=0}^{m} n_{i}$ lignes et $n_{\ell}$ colonnes
d\'efinie par
\begin{eqnarray*}
B_{\ell} & = & \left( b_{\ell;k,s} \right)_{\stackrel{0 \leq s \leq \sigma-1}{0 \leq k \leq n_{\ell}-1}}
= \left( \frac{1}{k!} \left( {{\rm d} \over {\rm d}X} \right)^{k} \left( X^{s} \right) (x_{\ell})
               \right)_{\stackrel{0 \leq s \leq \sigma-1}{0 \leq k \leq n_{\ell}-1}}.
\end{eqnarray*}
\end{defi}

On remarque que
$$
\frac{1}{k!} \left( {{\rm d} \over {\rm d}X} \right)^{k} \left( X^{s} \right) (x_{\ell})
= \left\{ \begin{array}{ll}
                    0  & \mbox{si $k > s$} \\
                    \displaystyle {s \choose k} x_{\ell}^{s-k} & \mbox{sinon.}
          \end{array}
  \right.
$$
Les d\'eterminants de Vandermonde g\'en\'eralis\'es vont jouer un r\^ole
important dans notre texte. Le lemme suivant permettra de montrer que le rang
de la matrice $\cM_{0}$ est \'egal \`a $S-1$.

\begin{lem}
\label{lem:gen-van}
Soit $\Delta$ un d\'eterminant de Vandermonde g\'en\'eralis\'e. On a 
$$
\Delta
= \prod_{0 \leq k < \ell \leq m} \left( x_{\ell}-x_{k} \right)^{n_{\ell}n_{k}}.
$$
\end{lem}

\begin{dem}
Voir \cite{El}.
\end{dem}

\subsection{Hauteurs de matrices}

Nous suivons ici le texte de F. Gramain \cite{Gr}. On y trouve plusieurs
identit\'es concernant certains d\'eterminants, utiles en th\'eorie des
nombres transcendants, ainsi que la notion de {\og hauteur de matrice\fg}
dont nous rappelons la d\'efinition.

\begin{defi}
Soient $m,n,d$ trois entiers strictement positifs. Soit $M$ une matrice
{\`a} $m$ lignes et $n$ colonnes, {\`a} coefficients dans le corps de nombres
$\bK$ de degr\'e $d$ sur $\bQ$, de rang $m \leq n$. On d\'efinit ses valeurs
absolues locales par:\\
$({\rm i})$ si $v$ est une place ultram\'etrique,
$$
|M|_{v}=\max_{I} \left| \Delta_{I} \right|_{v},
$$
o{\`u} le maximum porte sur toutes les parties $I$ {\`a} $m$ \'el\'ements de
$\{ 1, \ldots, n \}$, et o{\`u} $\Delta_{I}$ est le d\'eterminant $m\times{m}$
extrait de $M$ dont les colonnes sont les \'el\'ements de $I$,\\
$({\rm ii})$ si $v$ est une place archim\'edienne r\'eelle,
$$
|M|_{v}=|\det\left( M^{\ t}M \right)|_{v}^{1/2},
$$
o{\`u} $^{\ t}M$ est la matrice transpos\'ee de la matrice $M$,\\
$({\rm iii})$ si $v$ est une place archim\'edienne complexe,
$$
|M|_{v}=|\det\left(MM^{\star}\right)|_{v}^{1/2},
$$
o{\`u} $M^{\star}$ est la matrice adjointe de la matrice $M$.

La hauteur de la matrice $M$ est d\'efinie par:
$$
H(M) = \prod_{v \in M_\bK} |M|_{v}^{d_{v}/d}.
$$
\end{defi}

Pour une matrice carr\'ee r\'eguli\`ere, on a $H({}^{t}\!M) = H(M) = 1$, par la
formule du produit. Cela permet, pour une matrice $M$ dont le rang est \'egal
au nombre de colonnes, de d\'efinir
$$
H(M) = H({}^{t}\!M).
$$

\begin{lem}
\label{lem:3.9}
Pour toute place archim\'edienne, on a
$$
|M|_{v} = \left( \sum_{I} \left| \Delta_{I} \right|_{v}^{2} \right)^{1/2}.
$$
\end{lem}

\begin{dem}
Ceci provient essentiellement de la formule de Cauchy-Binet
(voir \cite{Gr} ou \cite{Se} par exemple).
\hfill$\Box$
\end{dem}

\vspace{3.0mm}

On d\'{e}signera par {\it matrice orthogonale $M^{\perp}$ \`{a}} $M$ toute
matrice \`{a} $n-m$ lignes et $n$ colonnes dont les lignes forment une
$\bK$-base de l'orthogonal du sous-espace de $\bK^{n}$ engendr\'{e} par les
lignes de $M$.

\begin{lem}
\label{lem:3.2}
Pour toute matrice $M^{\perp}$ orthogonale {\`a} $M$, on a 
$$
H\left(M^{\perp}\right)=H\left(M\right).
$$
\end{lem}

\begin{dem}
Voir \cite{StVa}.
\hfill$\Box$
\end{dem}

\vspace{3.0mm}

Le lien entre la hauteur de la matrice $\cM_{0}$ et celle de son orthogonal
sera utile; les composantes de cette matrice $\cM^{\perp}_{0}$ sont les
coefficients des polyn{\^o}mes de Hermite-Pad\'e. Cela permettra de ramener
la majoration de hauteurs de matrices \`a celle de hauteurs de polyn{\^o}mes.

\subsubsection{Estimations pr\'eliminaires de $H \left( \cM_{0} \right)$}

\begin{lem}
\label{lem:factorials}
{\rm (i)} Pour tout nombre r\'eel $x>1$, on a
$$
\sqrt{2 \pi (x-1)} \left( \frac{x-1}{e} \right)^{x-1} < \Gamma(x)
< \sqrt{2 \pi (x-1)} \left( \frac{x-1}{e} \right)^{x-1} e^{1/(12(x-1))}.
$$

{\rm (ii)} Pour tout entier naturel $n$, et tout entier $k$ dans l'intervalle
$0 \leq k\leq n$, on a
\begin{equation}
\label{eq:bin1}
{n \choose k} <  2^{n}\sqrt{ \frac{2}{\pi(n+1/2)} } 
\hspace{3.0mm} \mbox{ et } \hspace{3.0mm}
{n \choose k} \leq 2^{n}\sqrt{\frac{1}{n+1}}.
\end{equation}

Pour $n \geq 1$, et tout entier $k$ dans l'intervalle
$0 \leq k\leq n$, on a
\begin{equation}
\label{eq:bin2}
{n \choose k} \leq 2^{n}\sqrt{\frac{3}{4(n+1)}}.
\end{equation}
\end{lem}

\begin{dem}
(i) Cet encadrement de $\Gamma(x)$ r\'esulte de la formule de Stirling: voir
\cite[equation~6.1.38 on p.257]{AbSt}.

(ii) Pour d\'emontrer la premi\`ere in\'egalit\'e, on utilise la majoration du
th\'eor\`eme~1.1 de \cite{Ko} appliqu\'ee avec $\alpha=1/2$: pour tout entier
$n$ strictement positif, on a
$$
\frac{\Gamma(n+1/2)}{\Gamma(1/2) n!}
< \frac{1}{\Gamma(1/2)\sqrt{n+1/4}}.
$$

Pour $n$ pair, on a 
$$
\max_{k} {n \choose k} = {n \choose n/2} = 2^{n} \frac{\Gamma(n/2+1/2)}{\Gamma(1/2) (n/2)!}.
$$

L'in\'egalit\'e annonc\'ee r\'esulte de la relation $\Gamma(1/2)=\sqrt{\pi}$.

Pour $n$ impair, la valeur maximale de ${n \choose k}$ est atteinte pour
$k=(n+1)/2$, et on a
$$
2 {n \choose (n+1)/2} = {n+1 \choose (n+1)/2}
< 2^{n+1} \sqrt{\frac{2}{\pi(n+1+1/2)} }.
$$

Par cons\'equent l'in\'egalit\'e annonc\'ee est encore vraie pour $n$ impair.

La seconde in\'egalit\'e r\'esulte de la premi\`ere pour $n \geq 1$, et
un rapide calcul montre que cette estimation reste valable pour $n=0$.

Finalement, pour $n \geq (16-3\pi)/(6\pi-16)=2.307\ldots$, on a
$$
\sqrt{\frac{2}{\pi(n+1/2)}} \leq \sqrt{\frac{3}{4(n+1)}}.
$$
Donc la derni\`ere in\'egalit\'e r\'esulte de la premi\`ere pour $n \geq 3$, et
un calcul montre que cette derni\`ere estimation reste valable pour $n=1,2$.
\hfill$\Box$
\end{dem}

\begin{lem}
\label{lem:lambda-size}
{\rm (i)} Supposons que $m=L-1$ et que $n_{0} = \cdots = n_{L-1}=K$, o\`{u}
$K$ et $L$ d\'{e}signent des entiers avec $K \geq 1$ et $L \geq 2$.
Pour tout entier $0 \leq k \leq K-1$ et $0 \leq \ell \leq L-1$, le cardinal de
$\Lambda(\ell,k)$ est
$$
{K-k-1+L-2 \choose L-2}.
$$

{\rm (ii)} Pour $K \geq 1$ et $ L \geq 2$, on a 
$$
|\Lambda(\ell,k)| \leq 2^{K+L-3} \sqrt{\frac{2}{L}}.
$$
\end{lem}

\begin{dem}
(i) C'est Lemma~3.7.3 de \cite[p.33]{Ca}.

(ii) Pour $L=2$, on a $|\Lambda(k,\ell)|=1$, alors que
$$
2^{K+2-3}\sqrt{2/2} =2^{K-1} \geq 1
$$ 
pour $K \geq 1$.
Ainsi, pour terminer cette d\'emonstration, on peut supposer $L \geq 3$.

Comme, pour $y$ fix\'{e} positif, $\displaystyle{x \choose y}$ est croissante en
$x \geq 0$, et comme $0 \leq K-k-1+L-2 \leq K+L-3$, en appliquant la premi\`ere
in\'egalit\'e du lemme~\ref{lem:factorials} (ii), on obtient
$$
{K-k-1+L-2 \choose L-2} \leq {K+L-3 \choose L-2}
<2^{K+L-3} \sqrt{\frac{2}{\pi (K+L-3+1/2)}}.
$$

Pour $K \geq 1$ et $L>2$, on a $\pi (K+L-5/2) \geq \pi (L-3/2) \geq L$. Cela
compl\`ete la d\'emonstration du lemme~\ref{lem:lambda-size}.
\hfill$\Box$
\end{dem}

\begin{lem}
\label{lem:lem1}
Fixons deux entiers strictement positifs $b$ et $n$. Quand $a_{1}, \ldots, a_{n}$
et $A$ sont des nombres r\'eels positifs, on pose
$$
F \left( a_{1}, \ldots, a_{n} \right)
= \frac{1}{(b!)^{n}} \prod_{i=1}^{n} \prod_{j=1}^{b} \left( a_{i} + j \right)
$$
et
$$
g \left( a_{1}, \ldots, a_{n} \right) = a_{1} + \cdots + a_{n} - A.
$$
Alors, pour $A>0$, la valeur maximale de $F$ sur le domaine d\'efini par l'\'egalit\'e
$g \left( a_{1}, \ldots, a_{n} \right) =0$ est atteinte pour $a_{1}=\cdots=a_{n}=A/n$.
De plus, ce maximum est une fonction croissante de $A$.
\end{lem}

\begin{dem}
Posons
\begin{eqnarray*}
f \left( a_{1}, \ldots, a_{n} \right)
& = & \log \left( a_{1} + 1 \right) + \cdots + \log \left( a_{1} + b \right) \\
&   & + \cdots + \log \left( a_{n} + 1 \right) + \cdots + \log \left( a_{n} + b \right),
\end{eqnarray*}
de sorte que
\begin{displaymath}
f \left( a_{1}, \ldots, a_{n} \right)
= \log \left( (b!)^{n} F \left( a_{1}, \ldots, a_{n} \right) \right).
\end{displaymath}

Comme $b$ est fix\'e, il suffit de consid\'erer $f$. On applique la m\'ethode
des multiplicateurs de Lagrange pour trouver les valeurs extr\'emales de $f$.
On a 
\begin{displaymath}
\nabla f \left( a_{1}, \ldots, a_{n} \right)
= \left( \frac{1}{a_{1} + 1} + \cdots + \frac{1}{a_{1} + b},
\ldots,
\frac{1}{a_{n} + 1} + \cdots + \frac{1}{a_{n} + b} \right)
\end{displaymath}
et
\begin{displaymath}
\nabla g \left( a_{1}, \ldots, a_{n} \right)
= \left( 1, \ldots, 1 \right).
\end{displaymath}

La m\'ethode des multiplicateurs de Lagrange nous conduit \`{a}
r\'{e}soudre le syst\`{e}me d'\'{e}quations:
$$
\nabla f \left( a_{1}, \ldots, a_{n} \right)
= \lambda \nabla g \left( a_{1}, \ldots, a_{n} \right).
$$
avec la condition $g \left( a_{1}, \ldots, a_{n} \right)=0$. Ceci nous conduit
au syst\`{e}me d'\'{e}quations suivant:
\begin{eqnarray*}
\frac{1}{a_{1} + 1} + \cdots + \frac{1}{a_{1} + b} & = & \lambda \\
\cdots & & \\
\frac{1}{a_{n} + 1} + \cdots + \frac{1}{a_{n} + b} & = & \lambda \\
a_{1} + \cdots + a_{n} - A & = & 0,
\end{eqnarray*}
en les inconnues $a_{1}, \ldots, a_{n}, \lambda$.

La fonction $x\longmapsto\sum_{i=1}^{b} (x+i)^{-1}$ est strictement d\'{e}croissante
pour $x \geq 0$, donc les $n$ premi\`eres \'equations sont satisfaites si et seulement
si $a_{1}=\cdots=a_{n}$. Si $a$ est cette valeur commune, la derni\`ere \'equation
donne $na=A$, c'est-\`a-dire $a=A/n$. Le maximum 
$$
f(a,\ldots,a)=n\sum_{j=1}^b \log\left( \frac{A}{n}+j\right)
$$
est clairement une fonction croissante de $A$. D'o\`u le lemme~\ref{lem:lem1}.
\hfill
$\Box$
\end{dem}

\begin{lem}
\label{lem:lem2}
Soient $K$ et $L$ deux entiers strictement positifs. Pour tous $0 \leq k \leq K-1$
et $0 \leq \ell \leq L-1$, soit
$\gamma=\left( \gamma_{p}\right)_{\stackrel{0 \leq p \leq L-1}{p \neq \ell}}
$ un \'el\'ement de $\bZ^{L-1}$ satisfaisant 
$$
\displaystyle\sum_{\stackrel{p=0}{p \neq \ell}}^{L-1}\gamma_{p}=K-k-1
$$
et $\gamma_{p} \geq 0$ pour tout $0 \leq p \leq L-1$, $p \neq \ell$.
Alors
\begin{displaymath}
\prod_{p=0, p \neq \ell}^{L-1} { \gamma_{p} + K-1 \choose K-1}
\leq \left( e \min(K,L) \right)^{K-1}.
\end{displaymath}
\end{lem}

\begin{dem}
Si $K=1$, alors
\begin{displaymath}
\prod_{p=0, p \neq \ell}^{L-1} { \gamma_{p} + K-1 \choose K-1}
= \prod_{p=0, p \neq \ell}^{L-1} { \gamma_{p} \choose 0}
= 1 = (e\min(K,L))^{K-1},
\end{displaymath}
comme annonc\'e.

Si $L=1$, on a
\begin{displaymath}
\prod_{p=0, p \neq \ell}^{L-1} { \gamma_{p} + K-1 \choose K-1}
= 1 \leq (e\min(K,L))^{K-1},
\end{displaymath}
pour tout entier strictement positif $K$, comme voulu.

Pour la suite de cette d\'emonstration, on peut donc supposer que $K$ et $L$
sont $ \geq 2$. On applique le lemme~\ref{lem:lem1} avec $A=K-k-1$, $b=K-1$ et
$n=L-1$. Comme le maximum est une fonction croissante en $A$,
et comme $k$ est $ \geq 0$, on remplace $A$ par $K-1$ et on obtient 
\begin{eqnarray*}
\prod_{p=0, p \neq \ell}^{L-1} { \gamma_{p} + K-1 \choose K-1}
& \leq & \left( \frac{1}{(K-1)!} \prod_{i=1}^{K-1} \left( \frac{K-1}{L-1} + i \right) \right)^{L-1} \\
&   =  & \left( \frac{\displaystyle \Gamma \left( \frac{K-1}{L-1}+K \right)}
{\displaystyle \Gamma(K) \Gamma \left( \frac{K-1}{L-1}+1 \right)} \right)^{L-1}.
\end{eqnarray*}

On suppose dans un premier temps $L \leq K$. Les in\'egalit\'es $2 \leq L \leq K$
permettent d'appliquer le lemme~\ref{lem:factorials}(i):
$$
\Gamma \left( \frac{K-1}{L-1}+K \right)
< \sqrt{\frac{2\pi (K-1)L}{L-1}} \left( \frac{(K-1)L}{e(L-1)} \right)^{(K-1)L/(L-1)} e^{1/12},
$$
\begin{eqnarray*}
\sqrt{2\pi (K-1)} \left( \frac{K-1}{e} \right)^{(K-1)}
& < & \Gamma(K), \\
\sqrt{\frac{2\pi (K-1)}{L-1}} \left( \frac{K-1}{e(L-1)} \right)^{(K-1)/(L-1)}
& < & \Gamma \left( \frac{K-1}{L-1}+1 \right).
\end{eqnarray*}

Par cons\'equent
$$
\frac{\displaystyle \Gamma \left( \frac{K-1}{L-1}+K \right)}
{\displaystyle \Gamma(K) \Gamma \left( \frac{K-1}{L-1}+1 \right)}
< \sqrt{\frac{L}{2\pi (K-1)}} e^{1/12} \frac{L^{L(K-1)/(L-1)}}{(L-1)^{K-1}}.
$$

En utilisant les in\'egalit\'es 
$$
e^{1/12} \sqrt\frac{L}{2\pi (K-1)} \leq e^{1/12}\sqrt\frac{K}{2\pi (K-1)}<1
$$ 
et $(1+1/x)^{x}<e$ pour tout nombre r\'eel $x>0$, on obtient 
$$
\left( \frac{\displaystyle \Gamma \left( \frac{K-1}{L-1}+K \right)}
{\displaystyle \Gamma(K) \Gamma \left( \frac{K-1}{L-1}+1 \right)} \right)^{L-1}
< L^{K-1} \left[ \left(  1 + \frac{1}{L-1} \right)^{L-1} \right]^{K-1}
< (eL)^{K-1}.
$$
Ceci d\'emontre le lemme~\ref{lem:lem2} dans le cas $L \leq K$.

Supposons maintenant $L > K$. Dans ce cas, au moins $L-K$ nombres r\'eels parmi
les $\gamma_{p}$ doivent \^etre nuls. Ceux-ci peuvent \^etre ignor\'es et on
peut supposer $L=K$, cas pour lequel on vient de montrer que l'in\'egalit\'e
du lemme~\ref{lem:lem2} \'etait valable.
\hfill
$\Box$
\end{dem}

\subsubsection{Une majoration pour $H \left( \cM_{0} \right)$}

On reprend les notations de la premi\`ere section et on utilise
la matrice introduite \`a cette occasion.

Nous allons \'etablir une majoration pour la hauteur de la matrice $\cM_{0}$,
dont la pr\'ecision provient du lemme de dualit\'e~\ref{lem:3.2}. Cette estimation
jouera un r\^ole important dans la d\'emonstration du th\'eor\`eme principal.

\begin{lem}
\label{lem:m0-ortho}
Soient $K$ et $L$ des entiers strictement positifs avec $L \geq 2$. On peut
choisir pour matrice orthogonale \`{a} $\cM_{0}$ la matrice unilinge
$$
\cM_{0}^{\perp} = \left( \frac{p_{\ell,k}}{k!} \right)_{(k,\ell)},
$$
o\`{u} $(k,\ell)$ d\'{e}crit le carr\'{e} $0 \leq k \leq K-1$, $0 \leq \ell \leq L-1$.
\end{lem}

\begin{dem}
La matrice $\cM_{0}$ est de taille $(S-1) \times S$, o\`u $S=KL$.
La matrice obtenue en divisant par $k!$ chaque colonne $\cM_{0}$ index\'ee en
$(k,\ell)$, co\"{\i}ncide avec la matrice obtenue en rempla\c{c}ant la derni\`ere
ligne de la matrice de Vandermonde g\'en\'eralis\'ee avec $m=L-1$, $n_{\ell}=K$
et $x_{\ell}=\ell$ pour $\ell=0, \ldots, L-1$. Cette matrice a pour rang $S-1$:
en effet, la matrice de Vandermonde g\'en\'eralis\'ee a un d\'eterminant non nul,
elle est donc de rang $S$. Par cons\'equent $\cM_{0}^{\perp}$ est un vecteur.
Les composantes de ce vecteur $\cM_{0}^{\perp}$ sont les nombres $x_{\ell,k}$
tels que, pour tout $s=0, \ldots, S-2$,
\begin{equation}
\label{eq:herm}
\sum_{\ell=0}^{L-1} \sum_{k=0}^{K-1} x_{\ell,k}
\left( \frac{{\rm d}}{{\rm d}z} \right)^{s}
\left( z^{k}e^{\ell z} \right) \left( 0 \right)
= s! \sum_{\ell=0}^{L-1}
\sum_{k=0}^{\min(s,K-1)} x_{\ell,k} \frac{\ell^{s-k}}{(s-k)!}=0.
\end{equation}

Nous utiliserons les approximants de Hermite-Pad\'e pour de grandes valeurs de
$L$ avec $m=L-1$.

Pour $0\le \ell\le L-1$, consid\'erons le polyn\^ome 
$$
P_{\ell}(z)=\sum_{k=0}^{\min(s,K-1)} x_{\ell,k}z^{k}.
$$ 
On a 
$$
n_{\ell}=\deg(P_{\ell}(z))+1=\min(s+1,K).
$$ 
La somme
$$
\sum_{k=0}^{\min(s,K-1)}x_{\ell,k} {s \choose k} \ell^{s-k}
$$
est \'egale \`a $s!$ fois le coefficient de $z^{s}$ dans le d\'eveloppement de
$P_{\ell}(z) \cdot e^{\ell z}$ en $z=0$. D'apr\`es (\ref{eq:herm}), ce
coefficient est nul pour $s=0, \ldots, S-2$
et ainsi, ${\rm ord}_{z=0} \left( P_{\ell}(z) \cdot e^{\ell z} \right) \geq S-1$.
Notons que $\sigma=\sum_{\ell=0}^{m} n_{\ell} \leq KL=S$, d'o\`u
$S-1 \geq \sigma-1$.

Par cons\'equent les $x_{\ell,k}$ sont pr\'ecis\'ement les coefficients
$p_{\ell,k}/k!$ des approximants de Hermite-Pad\'e, pour les fonctions
$1,e^{z}, \ldots, e^{(L-1)z}$ et pour les param\`etres $n_{0} = \cdots = n_{L-1} = K$,
d\'etermin\'es lors de la section pr\'ec\'edente. On utilise le lemme~\ref{lem:2.2}
avec $m=L-1$, $n_{0}=\cdots=n_{L-1}=K$ et $x_{\ell}=\ell$ $(0 \leq \ell \leq L-1)$.
On trouve
$$
p_{\ell,k} = \sum_{\Lambda(\ell,k)}
\prod_{\stackrel{p=0}{p \neq \ell}}^{L-1}
\left( {(-1)^{\gamma_{p}} \over (\ell-p)^{\gamma_{p}+K}}
{\gamma_{p}+K-1 \choose K-1} \right),
$$
o{\`u} $\Lambda(\ell,k)$ a \'et\'e d\'efini dans la deuxi\`eme section.
\end{dem}

\begin{prop}
\label{prop:3.10}
Pour $K \geq 1$ et $L \geq 2$, on a 
$$
H \left( \cM_{0} \right)
\leq \frac{\sqrt{6}}{16L} 2^{KL+L} D_{K-1,L-1}
     \left( \frac{\sqrt{3}ed_{L-1}\min(K,L)}{2\sqrt{L}} \right)^{K-1}.
$$
\end{prop}

\begin{dem}
Pour $L \geq 2$, par les lemmes~\ref{lem:3.2} et \ref{lem:m0-ortho},
$$
H \left( \cM_{0} \right) = H \left( \cM_{0}^{\perp} \right)
= H \left( \left( {1 \over k!} p_{\ell,k} \right)_{k,\ell} \right).
$$ 

La hauteur d'une matrice ne change pas si on multiplie tous ses coefficients
par un m{\^e}me nombre non nul. On choisit le multiplicateur
$$
\rho = d_{L-1}^{K-1}(L-1)!^{K} D_{K-1,L-1}
$$
et, pour tout $0 \leq k \leq K-1$ et $0 \leq \ell \leq L-1$, on consid\`{e}re
les nombres
$$
y_{\ell, k}
= \rho \frac{p_{\ell,k}}{k!}
= \pm \frac{D_{K-1,L-1}}{k!} d_{L-1}^{K-1} {L-1 \choose \ell}^{K}
\sum_{\Lambda(\ell,k)}
\prod_{\stackrel{p=0}{p \neq \ell}}^{L-1}
\left( \frac{(-1)^{\gamma_{p}}}{(\ell-p)^{\gamma_{p}}}
{\gamma_{p}+K-1 \choose K-1} \right).
$$

On observe que pour tout nombre premier $q \leq L-1$, on a $q^{k}| d_{L-1}^{k}$
et $v_{q}(k!) \leq k$. On en d\'eduit $D_{K-1,L-1} d_{L-1}^{k}/k! \in \bZ$. On
utilise le fait que ${\gamma_{p}+K-1 \choose K-1} \in \bZ$.
La condition $\sum_{p \neq \ell} \gamma_{p} = K-k-1$ implique
$$
d_{L-1}^{K-k-1}
\prod_{\stackrel{p=0}{p \neq \ell}}^{L-1}
\frac{1}{(\ell-p)^{\gamma_{p}}} \in \bZ.
$$

D'o\`u $y_{\ell,k} \in \bZ$. On obtient 
\begin{eqnarray*}
& & H \left( \left( {1 \over k!} p_{\ell, k} \right)_{k,\ell} \right)
= H \left( \left( y_{\ell, k} \right)_{k,\ell} \right)
= \prod_{v \in M_\bQ} \max_{k,\ell} \left| y_{\ell, k} \right|_{v} \\
& = & \max_{k,\ell} \left| y_{\ell, k} \right|
\prod_{p \in \cP} \max_{k,\ell} \left| y_{\ell, k} \right|_{p}
\leq \max_{k,\ell} \left| y_{\ell, k} \right| \\
& \leq & \max_{k,\ell} \left| \frac{D_{K-1,L-1}}{k!} d_{L-1}^{K-1} {L-1 \choose \ell}^{K}
\sum_{\Lambda(\ell,k)}
\prod_{\stackrel{p=0}{p \neq \ell}}^{L-1}
\left( {(-1)^{\gamma_{p}} \over (\ell-p)^{\gamma_{p}}}
{\gamma_{p}+K-1 \choose K-1} \right) \right|.
\end{eqnarray*}

Des lemmes~\ref{lem:lem2}, \ref{lem:lambda-size}(ii) et la derni\`{e}re
in\'{e}galit\'{e} du \ref{lem:factorials}(ii), on d\'eduit
\begin{eqnarray*}
& &
\left| \frac{D_{K-1,L-1}}{k!} d_{L-1}^{K-1} {L-1 \choose \ell}^{K}
\sum_{\Lambda(\ell,k)}
\prod_{\stackrel{p=0}{p \neq \ell}}^{L-1}
\left( {(-1)^{\gamma_{p}} \over (\ell-p)^{\gamma_{p}}}
{\gamma_{p}+K-1 \choose K-1} \right) \right| \\
& \leq & D_{K-1,L-1} d_{L-1}^{K-1} {L-1 \choose \ell}^{K} |\Lambda(\ell,k)|
         (e\min(K,L))^{K-1} \\
& \leq & D_{K-1,L-1} \left( e\min(K,L) d_{L-1} {L-1 \choose \ell} \right)^{K-1}
         {L-1 \choose \ell} \sqrt{\frac{2}{L}} 2^{K+L-3} \\
& <    & D_{K-1,L-1} \left( \frac{\sqrt{3}}{2\sqrt{L}}e\min(K,L) d_{L-1} \right)^{K-1}
         2^{K(L-1)} \frac{\sqrt{6}}{2L} 2^{K+L-3} \\
& <    & D_{K-1,L-1} \left( \frac{\sqrt{3}}{2\sqrt{L}}e\min(K,L) d_{L-1} \right)^{K-1}
         2^{KL+L} \frac{\sqrt{6}}{16L},
\end{eqnarray*}
ce qui compl\`ete la d\'emonstration de la proposition~\ref{prop:3.10}.
\hfill$\Box$
\end{dem}

\begin{rem}
Pour $K \geq 1$ et $L=1$, on a $\cM_{0}^{\perp}=(0,\ldots,0,1)$  et $H(\cM_{0})=1$.
\end{rem}

\subsection{Lemme de Schwarz}

Le lemme ci-dessous fournit une majoration analytique du d\'eterminant \'etudi\'e:

\begin{lem}[Lemme de Schwarz]
\label{lem:Sc}
Soient $T$ un entier positif, $r$ et $R$ deux nombres r\'eels v\'erifiant
$0 < r \leq R$ et $\psi$ une fonction d'une variable complexe analytique
dans le disque $|z| \leq R$. Supposons que $\psi$ a un z\'ero de multiplicit\'e
au moins $T$ en $0$. Alors
$$
|\psi|_{r} \leq \left( \frac{R}{r} \right)^{-T} |\psi|_{R},
$$
o\`u $|\psi|_{R}=\max_{|z|=R}|\psi(z)|$.
\end{lem}

\begin{dem}
Voir \cite{Wa}, paragraphe~2.2.3, lemme~2.4, page~37.
\hfill$\Box$
\end{dem}

\vspace{3.0mm}

Voici une majoration du d\'eterminant $\cD$, introduit dans la
d\'efinition~\ref{defi:bigF}. On utilise les nombres $w_{k,\ell}$ introduits
dans la d\'efinition~\ref{defi:wkl}. 

\begin{lem}
\label{lem:4.10}
Soient $K$ et $L$ deux entiers strictement positifs avec $L \geq 2$, $\mu$ un
entier $\geq 0$ et $E$ un nombre r\'eel $\geq 1$. Soit $\cN$ un nombre
r\'eel strictement positif tel que 
\begin{equation}
\label{eq:lem4.10a}
\max\left\{
\sum_{k=0}^{K-1}\sum_{\ell=0}^{L-1} \left| w_{k,\ell} \right|, \hspace{1.0mm}
\sum_{k=0}^{K-1}\sum_{\ell=0}^{L-1} \left| \Phi_{k,\ell}^{(\mu)}(z\beta) \right|_{E}
\right\}
\leq e^{\cN}.
\end{equation}

Supposons que le nombre $\epsilon= \left| \exp(\beta)-\alpha \right|$ v\'erifie
\begin{equation}
\label{eq:lem4.10b}
\epsilon < E^{-KL}.
\end{equation}

Alors
$$
\log \left| \cD \right|
\leq -(KL-\mu-1)\log E + \cN + \log \left| \cM_{0} \right| + \log(2).
$$
\end{lem}

\begin{dem}
Consid\'erons la fonction
$$
\cD(z) = \det \left(
\begin{array}{l}
    \cM_{0} \\
    \left( \Phi_{k,\ell}^{(\mu)}(\beta z) + \epsilon w_{k,\ell}
    \right)_{0 \leq k \leq K-1,0 \leq \ell \leq L-1}
\end{array} \right),
$$
de sorte que
$$
\cD = \cD(1).
$$

Alors
\begin{equation}
\label{eq:*}
\cD(z) = \cD_{0}(z) + \epsilon \cD_{1}(z),
\end{equation}
o{\`u}
$$
\cD_{0}(z) = \det \left(
\begin{array}{l}
    \cM_{0} \\
    \left( \Phi_{k,\ell}^{(\mu)}(\beta z) \right)_{0 \leq k \leq K-1,0 \leq \ell \leq L-1}
\end{array} \right)
$$
et
$$
\cD_{1}(z) = \det \left(
\begin{array}{l}
    \cM_{0} \\
    \left( w_{k,\ell} \right)_{0 \leq k \leq K-1,0 \leq \ell \leq L-1}
\end{array}\right).
$$

On a 
$$
\left( \frac{\partial}{\partial z} \right)^{t} \cD_{0}(z) = \det \left(
\begin{array}{l}
    \cM_{0} \\
    \left( \beta^{t} \Phi_{k,\ell}^{(\mu+t)}(\beta z) \right)_{0 \leq k \leq K-1,0 \leq \ell \leq L-1}
\end{array} \right).
$$

Supposons que $t+\mu \leq S-2=KL-2$. Alors la ligne 
$$
\left( \beta^{t} \Phi_{k,\ell}^{(\mu+t)}(0 \cdot \beta) \right)_{0 \leq k \leq K-1,0 \leq \ell \leq L-1}
$$
est \'egale \`a $\beta^{t}$ fois la ligne dans $\cM_{0}$ index\'ee par $t+\mu$.
Donc, pour $0 \leq t< KL-1-\mu$,
$$
\left( \frac{\partial}{\partial z} \right)^{t} \cD_{0}(0) = 0.
$$

Par cons\'equent le d\'eterminant $\cD_{0}(z)$ a un z\'ero \`a l'origine d'ordre
sup\'erieur ou \'egal \`a $KL-\mu-1$.

Du Lemme de Schwarz (Lemme~\ref{lem:Sc}),
on d\'eduit, pour tout $E \geq 1$,
\begin{equation}
\label{eq:**}
\log \left| \cD_{0}(1) \right|
\leq -(KL-\mu-1) \log E + \log \left| \cD_{0}(z) \right|_{E}.
\end{equation}

Comme $\cD_{1}(z)$ ne d\'epend pas de $z$, on peut \'ecrire
$$
\log \left| \cD_{1}(1) \right| = \log \left| \cD_{1}(z) \right|_{E}.
$$

\'Etudions maintenant, pour $I=0$ et $I=1$, le terme
$\log \left| \cD_{I}(z) \right|_{E}$.
En d\'eveloppant $\cD_{I}(z)$ suivant la derni\`ere ligne, on obtient
\begin{eqnarray*}
\left| \cD_{I}(z) \right|
& \leq & \left\{
\begin{array}{ll}
    \displaystyle \sum_{\stackrel{0 \leq k \leq K-1}{0 \leq \ell \leq L-1}}
    \left| \Phi_{k,\ell}^{(\mu)}(\beta z) \right| \left| \Delta_{k,\ell} \right| &
    {\rm si} \hspace{3.0mm} I=0 \\
    \hspace{3.0mm} \\
    \displaystyle \sum_{\stackrel{0 \leq k \leq K-1}{0 \leq \ell \leq L-1}}
    \left| w_{k,\ell} \right| \left| \Delta_{k,\ell} \right| &
    {\rm si} \hspace{3.0mm} I=1,
\end{array} \right.
\end{eqnarray*}
o\`u $\left| \Delta_{k,\ell} \right|$ est le mineur de la matrice $\cM_{0}$
dont la colonne d'indice $(k,\ell)$ a \'et\'e {\^o}t\'ee.

De l'in\'egalit\'e de Cauchy-Schwarz, on d\'eduit
\begin{eqnarray*}
\left| \cD_{I}(z) \right|_{E}^{2}
& \leq &
\max \left( \sum_{\stackrel{0 \leq k \leq K-1}{0 \leq \ell \leq L-1}}
\left| w_{k,\ell} \right|^{2},
\sum_{\stackrel{0 \leq k \leq K-1}{0 \leq \ell \leq L-1}} \left| \Phi_{k,\ell}^{(\mu)}(\beta z) \right|_{E}^{2} \right)
\left( \sum_{\stackrel{0 \leq k \leq K-1}{0 \leq \ell \leq L-1}} \left| \Delta_{k,\ell} \right|^{2} \right) \\
& \leq &
\max \left( \sum_{\stackrel{0 \leq k \leq K-1}{0 \leq \ell \leq L-1}}
\left| w_{k,\ell} \right|,
\sum_{\stackrel{0 \leq k \leq K-1}{0 \leq \ell \leq L-1}} \left| \Phi_{k,\ell}^{(\mu)}(\beta z) \right|_{E} \right)^{2}
\left( \sum_{\stackrel{0 \leq k \leq K-1}{0 \leq \ell \leq L-1}} \left| \Delta_{k,\ell} \right|^{2} \right).
\end{eqnarray*}

De plus, le lemme~\ref{lem:3.9} donne
$$
\left| \cM_{0} \right|^{2}
= \sum_{\stackrel{0 \leq k \leq K-1}{0 \leq \ell \leq L-1}} \left| \Delta_{k,\ell} \right|^{2}.
$$

On en d\'eduit, en utilisant l'hypoth\`ese (\ref{eq:lem4.10a}),
$$
\log \left| \cD_{I}(z) \right|_{E}
\leq \log \left| \cM_{0} \right| + \cN.
$$

En utilisant l'in\'egalit\'e (\ref{eq:**}), on obtient
$$
\log \left| \cD_{I}(1) \right| \leq \left\{
\begin{array}{ll}
    -(KL-\mu-1) \log E + \cN + \log \left| \cM_{0} \right|
    & {\rm si} \hspace{1.0mm} I=0, \\
    & \\
    \cN + \log \left| \cM_{0} \right|
    & {\rm si} \hspace{1.0mm} I=1.
\end{array}
\right.
$$

De l'\'egalit\'e (\ref{eq:*}), on d\'eduit
$$
\left| \cD \right| = \left| \cD(1) \right|
\leq 2 \max \left\{ \left| \cD_{0}(1) \right|,
                    \epsilon \left| \cD_{1}(1) \right|
            \right\}.
$$

Comme $E \geq 1$, en utilisant l'in\'egalit\'e (\ref{eq:lem4.10b}) on trouve
$$
\log \epsilon \leq -KL\log E \leq -(KL-\mu-1)\log E.
$$ 

Par cons\'equent
$$
\log \left| \cD \right|
\leq -(KL-\mu-1)\log E + \cN + \log \left| \cM_{0} \right| + \log(2).
$$

Ceci d\'emontre le lemme~\ref{lem:4.10}.
\hfill$\Box$
\end{dem}

\section{D\'emonstration du th\'eor\`eme~\ref{theo:4.6bisalg}}

Nous devons traiter le cas $\alpha=0$ s\'epar\'ement: dans ce cas, nous ne
pouvons pas utiliser notre lemme de z\'eros (proposition~\ref{prop:NW}). Nous
voulons minorer $|e^{\beta}|$. Notons que
$|e^{\beta}|=|e^{\Re{e}(\beta)}|>e^{-|\beta|}$. De (\ref{eq:theoTalg}), on en
d\'eduit que $KL \log (E) > |\beta|$, donc que $|e^{\beta}|>E^{-KL}$. Ainsi 
le th\'eor\`eme~\ref{theo:4.6bisalg} est vrai pour $\alpha=0$.

Nous supposerons dans la suite $\alpha \neq 0$. Commen\c{c}ons par d\'efinir
une nouvelle quantit\'e, $\cG_{\beta,\alpha}$, li\'ee \`a $\cD$.

\subsection{D\'{e}finition de $\cG_{\beta,\alpha}$}

Avec $\mu$ tel que d\'{e}fini dans D\'{e}finition~\ref{defi:mu}, on introduit
l'op\'erateur de d\'erivation
$$
\cF_{\mu}(\delta) = \frac{\delta(\delta-1) \cdots (\delta-\mu+1)}{\mu!},
$$
avec $\cF_{0}(\delta)=1$.

Consid\'erons le polyn\^ome en deux variables $\cG_{1}(X,Y)$ d\'efini par
\begin{equation}
\label{eq:defG1}
\cG_{1}(X,Y)
= d_{\mu}^{K-1} \left( \prod_{p \in \cP} \left| \cM_{0} \right|_{p} \right)
\cF_{\mu}(\delta) \left( \cH(X,Y) \right),
\end{equation}
o\`{u} $\cH(X,Y)$ est donn\'e par (\ref{eq:bigH}) et on a pos\'{e} par
convention $d_{\mu}=1$ quand $\mu=0$.

\begin{lem}
\label{lem:comp}
Le polyn\^ome $\cG_{1}(X,Y)$ est \`a coefficients dans $\bZ$. Alors
$$
\cG_{1}(\beta,\alpha) = {d_{\mu}^{K-1} \over \mu!}
\left( \prod_{p \in \cP} \left| \cM_{0} \right|_{p} \right) F(\beta,\alpha),
$$
o\`{u} $F(X,Y)$ est tel que d\'{e}fini dans d\'{e}finition~$\ref{defi:bigF}$,
et $\cG_{1}(\beta, \alpha) \neq 0$.
\end{lem}

\begin{dem}
Avec $0^0=1$, pour tout couple d'indices $(k,\ell)$, on a
$$
\delta^{i}(X^{k}Y^{\ell})
= \sum_{j=0}^{\min(i,k)} {i \choose j} \frac{k!}{(k-j)!} \ell^{i-j} X^{k-j}Y^{\ell},
$$
$$
\delta^{i}(X^{k})
= \left\{
\begin{array}{ll}
\displaystyle \frac{k!}{(k-i)!} X^{k-i} & \mbox{si $i \leq k$} \\
0 & \mbox{sinon,}
\end{array}
\right.
$$
et
$$
\cF^{(j)}_{\mu}(\ell) = \sum_{i=j}^{\mu} \lambda_{i,\mu} {i! \over (i-j)!} \ell^{i-j}
 = \sum_{i=j}^{\mu} \lambda_{i,\mu} (i-j+1) \cdots i \ell^{i-j}.
$$

Alors
\begin{eqnarray*}
\cF_{\mu}(\delta)(X^{k}Y^{\ell})
& = & \sum_{i=0}^{\mu} \lambda_{i,\mu}
          \sum_{j=0}^{\min(i,k)} {i \choose j} {k! \over (k-j)!} \ell^{i-j}X^{k-j}Y^{\ell} \\
& = & \sum_{i=0}^{\mu} \lambda_{i,\mu}
          \sum_{j=0}^{\min(i,k)} {i! \over (i-j)!} {k \choose j} \ell^{i-j}X^{k-j}Y^{\ell} \\
& = & \sum_{i=0}^{\mu} \lambda_{i,\mu}
          \sum_{j=0}^{\min(i,k)} (i-j+1) \cdots i {k \choose j} \ell^{i-j}X^{k-j}Y^{\ell} \\
& = & \sum_{i=0}^{\mu} \lambda_{i,\mu}
          \sum_{j=0}^{k} (i-j+1) \cdots i {k \choose j} \ell^{i-j}X^{k-j}Y^{\ell} \\
& = & \sum_{j=0}^{k} {k \choose j}
           \sum_{i=0}^{\mu} \lambda_{i,\mu} (i-j+1) \cdots i \ell^{i-j} X^{k-j}Y^{\ell} \\
& = & Y^{\ell} \sum_{j=0}^{k} {k \choose j}
           \left (\sum_{i=0}^{\mu} \lambda_{i,\mu} (i-j+1) \cdots i \ell^{i-j} \right) X^{k-j} \\
& = & Y^{\ell} \sum_{j=0}^{k} {k \choose j} \cF^{(j)}_{\mu}(\ell) X^{k-j}.
\end{eqnarray*}

Gr\^ace au lemme~\ref{lem:fel3}, on en d\'eduit que
$d_{\mu}^{K-1}\cF_{\mu}(\delta)(X^{k}Y^{\ell})$ est un polyn\^ome \`a coefficients
dans $\bZ$.

Comme
$$
\cF_{\mu}(\delta) \left( \cH(X,Y) \right)
=\det \left(
\begin{array}{l}
    \cM_{0} \\
    \left( \cF_{\mu}(\delta) \left( X^{k}Y^{\ell} \right)
    \right)_{0 \leq k \leq K-1,0 \leq \ell \leq L-1}
\end{array} \right),
$$
il r\'esulte de la d\'efinition de $\left| \cM_{0} \right|_{p}$ que les
coefficients du polyn\^ome $\cG(X,Y)$ sont dans $\bZ$.

En utilisant d'une part la d\'efinition de $\cH$ et d'autre part le fait
que le polyn\^ome en deux variables
$$
\det \left(
\begin{array}{l}
\cM_{0} \\
\left( \delta^{j} \left( X^{k}Y^{\ell} \right) \right)_{0 \leq k \leq K-1,0 \leq \ell \leq L-1}
\end{array} \right)
$$
s'annule au point $(\beta,\alpha)$ pour $0 \leq j < \mu$, on d\'eduit
\begin{eqnarray*}
\cF_{\mu}(\delta) \left( \cH \right)(\beta,\alpha)
& = & \det \left(
\begin{array}{l}
    \cM_{0} \\
    \left( \cF_{\mu}(\delta) \left( X^{k}Y^{\ell} \right)
    \right)_{0 \leq k \leq K-1,0 \leq \ell \leq L-1}
\end{array} \right) (\beta,\alpha) \\
& = & \det \left(
\begin{array}{l}
    \cM_{0} \\
    \left( \displaystyle{1 \over \mu!} \delta^{\mu} \left( X^{k}Y^{\ell} \right)
    \right)_{0 \leq k \leq K-1,0 \leq \ell \leq L-1}
\end{array} \right)(\beta,\alpha) \\
& & + \sum_{j=0}^{\mu-1} \lambda_{j,\mu}
\det \left(
\begin{array}{l}
\cM_{0} \\
\left( \delta^{j} \left( X^{k}Y^{\ell} \right) \right)_{0 \leq k \leq K-1,0 \leq \ell \leq L-1}
\end{array} \right)(\beta,\alpha) \\
& = & \det \left(
\begin{array}{l}
    \cM_{0} \\
    \left( \displaystyle{1 \over \mu!} \delta^{\mu} \left( X^{k}Y^{\ell} \right)
    \right)_{0 \leq k \leq K-1,0 \leq \ell \leq L-1}
\end{array} \right)(\beta,\alpha) \\
& = & {\delta^{\mu} \over \mu!} \cH(\beta,\alpha) = {1 \over \mu!} F(\beta,\alpha).
\end{eqnarray*}

Ainsi, les deux polyn\^omes $\displaystyle{1 \over \mu!} F(X,Y)$ et
$\cF_{\mu} (\delta) \left( \cH(X,Y) \right)$ prennent la m\^eme valeur au point
$(\beta,\alpha)$. Par la d\'{e}finition~\ref{defi:mu},
$\delta^{\mu}\cH(\beta, \alpha) \neq 0$ et la preuve du lemme~\ref{lem:comp}
est ainsi compl\`ete.
\hfill$\Box$
\end{dem}

\vspace{3.0mm}

On d\'efinit
\begin{equation}
\label{eq:defG2}
\cG_{2}(X,Y)
= \left( \prod_{p \in \cP} \left| \cM_{0} \right|_{p} \right)
\delta^{\mu} \left( \cH(X,Y) \right).
\end{equation}

Puisque $F(X,Y)=\delta^{\mu} \left( \cH(X,Y) \right)$, on a $\cG_{2}(X,Y) \in \bZ[X,Y]$
et
$$
\cG_{2}(\beta,\alpha)=\left( \prod_{p \in \cP} \left| \cM_{0} \right|_{p} \right)
F(\beta,\alpha).
$$

Finalement, on d\'efinit
\begin{equation}
\label{eq:defG}
\cG_{\beta,\alpha} = \min \left( \left| \cG_{1}(\beta,\alpha) \right|,
\left| \cG_{2}(\beta,\alpha) \right| \right)
= \min \left( 1, \frac{d_{\mu}^{K-1}}{\mu!} \right)
\left( \prod_{p \in \cP} \left| \cM_{0} \right|_{p} \right) \left| F(\beta,\alpha) \right|.
\end{equation}

\subsection{Minoration de $\cG_{\beta,\alpha}$}

\begin{defi}
Soit $\displaystyle f(X_{1},\ldots,X_{n})
=\sum_{i_{1},\ldots,i_{n}}a_{i_{1},\ldots,i_{n}}X_{1}^{i_{1}} \cdots X_{n}^{i_{n}}$
un polyn\^ome \`a coefficients complexes. On d\'efinit la longueur $\length(f)$
de $f$ comme la somme des valeurs absolues de ses coefficients, \`a savoir
$$
\length(f)=\sum_{i_{1},\ldots,i_{n}}|a_{i_{1},\ldots,i_{n}}|.
$$
\end{defi}

Pour minorer $\cG_{\beta,\alpha}$, on utilisera l'estimation suivante
de la longueur de $\cG_{1}$et $\cG_{2}$.

\begin{lem}
\label{lem:4.3.2}
$({\rm i})$ La longueur du polyn\^ome $\cG_{1}$ est major\'ee par
$$
\length \left( \cG_{1} \right) \leq d_{\mu}^{K-1} H \left( \cM_{0} \right) 2^{\mu+K-1} e^{L-1} \sqrt{L}.
$$

$({\rm ii})$La longueur du polyn\^ome $\cG_{2}$  est major\'ee par
$$
\length \left( \cG_{2} \right) \leq \mu! H \left( \cM_{0} \right) 2^{\mu+K-1} e^{L-1} \sqrt{L}.
$$
\end{lem}

\begin{dem}
(i) On commence par majorer $\length\left(\cF_{\mu}(\delta)(\cH(X,Y))\right)$.
Soit $\cM_{0,k,\ell}$ la matrice obtenue \`a partir de $\cM_{0}$ en supprimant
la colonne d'indice $(k,\ell)$. En utilisant la majoration $\length(f \pm g)
\leq \length(f)+\length(g)$ et en d\'eveloppant le d\'eterminant
$\cF_{\mu}(\delta)(\cH(X,Y))$ suivant la derni\`ere ligne, on obtient
\begin{eqnarray*}
\length \left( \cF_{\mu}(\delta)(\cH(X,Y)) \right)
& \leq & \sum_{k=0}^{K-1} \sum_{\ell=0}^{L-1}
         \length \left( \cF_{\mu}(\delta) \left( X^{k}Y^{\ell} \right) \det \left( \cM_{0,k,\ell} \right) \right) \\
& = & \sum_{k=0}^{K-1} \sum_{\ell=0}^{L-1}
      \length \left( \left( \sum_{j=0}^{\mu} \lambda_{j,\mu}
      \delta^{j} \left( X^{k}Y^{\ell} \right) \right)
      \det \left( \cM_{0,k,\ell} \right) \right) \\
& = & \sum_{k=0}^{K-1} \sum_{\ell=0}^{L-1}
      \length \left( \left( \sum_{j=0}^{\mu} \lambda_{j,\mu}
      \sum_{h=0}^{\min(j,k)} {j \choose h}{k! \over (k-h)!} \ell^{j-h} X^{k-h} Y^{\ell}
      \right) \right. \\
&   & \left. \hspace{7.0mm} \det \left( \cM_{0,k,\ell} \right) \displaystyle \right) \\
& = & \sum_{k=0}^{K-1} \sum_{\ell=0}^{L-1} \left| \det \left( \cM_{0,k,\ell} \right) \right|
           \cdot \sum_{j=0}^{\mu} \sum_{h=0}^{\min(j,k)}
           \left| \lambda_{j,\mu} \right| {j \choose h} {k! \over (k-h)!} \ell^{j-h}.
\end{eqnarray*}

En utilisant, de plus, l'in\'egalit\'e de Cauchy-Schwarz, on obtient
\begin{eqnarray*}
\length \left( \cF_{\mu}(\delta) (\cH(X,Y)) \right)
& \leq & \left( \sum_{k=0}^{K-1} \sum_{\ell=0}^{L-1}
               \left| \det \left( \cM_{0,k,\ell} \right) \right|^{2} \right)^{1/2} \\
& & \cdot \left( \sum_{k=0}^{K-1} \sum_{\ell=0}^{L-1}
       \left( \sum_{j=0}^{\mu} \sum_{h=0}^{\min(j,k)}
       \left| \lambda_{j,\mu} \right| {j \choose h} \frac{k!}{(k-h)!}
       \ell^{j-h} \right)^{2} \right)^{1/2}.
\end{eqnarray*}

On majore $\displaystyle{K-1 \choose h}$ par $\displaystyle \frac{2^{K-1}}{\sqrt{K}}$,
gr\^ace la seconde in\'egalit\'e du lemme~\ref{lem:factorials}(ii),
et $\displaystyle\sum_{h=0}^{\min(j,k)} {(L-1)^{j-h} \over (j-h)!}$ par
$e^{L-1}$, puis on utilise le lemme~\ref{lem:fel2} et le lemme~\ref{lem:3.9}
afin de conclure de la fa\c{c}on suivante:
\begin{eqnarray*}
\length \left( \cF_{\mu}(\delta) (\cH(X,Y)) \right)
& \leq & \left| \cM_{0} \right|
               \sqrt{KL} \sum_{j=0}^{\mu} \sum_{h=0}^{\min(j,k)}
               \left| \lambda_{j,\mu} \right| j! {K-1 \choose h} {(L-1)^{j-h} \over (j-h)!} \\
& \leq & \left| \cM_{0} \right| \sqrt{L} \sum_{j=0}^{\mu} \left| \lambda_{j,\mu} \right|
               j! 2^{K-1} e^{L-1} \\
& \leq & \left| \cM_{0} \right| 2^{\mu+K-1} e^{L-1} \sqrt{L}.
\end{eqnarray*}

La d\'efinition~(\ref{eq:defG1}) de $\cG_{1}(X,Y)$ termine la d\'emonstration de (i).

(ii) On proc\'ede de fa\c{c}on similaire: on d\'eveloppe tout d'abord le d\'eterminant
$\delta^{\mu}(\cH(X,Y))$ suivant la derni\`ere ligne, puis on d\'eveloppe
$\delta^{\mu} \left( X^{k}Y^{\ell} \right)$ et enfin on applique l'in\'egalit\'e
de Cauchy-Schwarz pour obtenir
\begin{eqnarray*}
\length \left( \delta^{\mu} \left( \cH(X,Y) \right) \right)
& \leq & \left( \sum_{k=0}^{K-1} \sum_{\ell=0}^{L-1}
               \left| \det \left( \cM_{0,k,\ell} \right) \right|^{2} \right)^{1/2} \\
&      & \cdot \left( \sum_{k=0}^{K-1} \sum_{\ell=0}^{L-1}
         \left( \sum_{h=0}^{\min(\mu,k)} {\mu \choose h} \frac{k!}{(k-h)!}
         \ell^{\mu-h} \right)^{2} \right)^{1/2}.
\end{eqnarray*}

En utilisant la majoration $\displaystyle{K-1 \choose h} \leq \displaystyle \frac{2^{K-1}}{\sqrt{K}}$,
$\displaystyle\sum_{h=0}^{\min(\mu,k)} \frac{(L-1)^{\mu-h}}{(\mu-h)!} \leq
e^{L-1}$ et le lemme~\ref{lem:3.9}, on obtient
\begin{eqnarray*}
\length \left( \delta^{\mu} (\cH(X,Y)) \right)
& \leq & \left| \cM_{0} \right|
               \left( KL \left( \sum_{h=0}^{\min(\mu,k)}
               \mu! {K-1 \choose h} {(L-1)^{\mu-h} \over (\mu-h)!} \right)^{2} \right)^{1/2} \\
& \leq & \left| \cM_{0} \right| \sqrt{L} \mu! 2^{K-1} e^{L-1}.
\end{eqnarray*}

On multiplie \'egalement le membre de droite de cette expression par $2^{\mu}$.
La d\'efinition~(\ref{eq:defG2}) de $\cG_{2}(X,Y)$ permet de compl\'eter cette
d\'emonstration du lemme~\ref{lem:4.3.2}.
\hfill$\Box$
\end{dem}

\begin{lem}[Liouville]
\label{lem:4.3.3}
Soient $\alpha_{1}, \ldots, \alpha_{n}$ des nombres alg\'ebriques. On pose 
$D=[\bQ(\alpha_{1}, \ldots, \alpha_{n}):\bQ]/[\bR(\alpha_{1}, \ldots, \alpha_{n}):\bR]$.
Soit $f$ un polyn\^ome de $\bZ[X_{1},\ldots,X_{n}]$, de degr\'e au plus
$N_{i}$ par rapport \`a la variable $X_{i}$, et qui ne s'annule pas au point
$(\alpha_{1}, \ldots, \alpha_{n})$. Alors
$$
\log |f(\alpha_{1},\ldots,\alpha_{n})|
\geq -(D-1)\log \length(f)+\sum_{i=1}^{n}N_{i} \log (\max(1, |\alpha_{i}|))
-D\sum_{i=1}^{n}N_{i}\hgt(\alpha_{i}).
$$
\end{lem}

\begin{dem}
Ceci est la g\'en\'eralisation de la version de l'in\'egalit\'e de Liouville
utilis\'ee page~298 de \cite{LMN} appliqu\'ee avec l'in\'egalit\'e
$|f| \leq \length(f)$. Voir aussi \cite{FeNe} Chap.~1 \S~1.7 et \cite{Wa}
Prop.~3.14. 
\hfill$\Box$
\end{dem}

\vspace{3.0mm}

Gr\^ace \`a la proposition~\ref{prop:3.10}, on en d\'eduit la minoration
suivante de $\cG_{\beta,\alpha}$:

\begin{prop}
\label{prop:trivalg}
Si $\alpha \neq 0$, on a:
\begin{eqnarray*}
\log \cG_{\beta,\alpha}
& \geq & -(D-1) \log \left( \min \left( d_{\mu}^{K-1}, \mu! \right) H \left( \cM_{0} \right)
                2^{\mu+K-1} e^{L-1}\sqrt{L} \right) \\
&      & -(K-1) \log (\cB) - (L-1)\log (\cA).
\end{eqnarray*}
\end{prop}

\begin{dem}
Gr\^ace aux lemmes~\ref{lem:4.3.2} et~\ref{lem:4.3.3},
en utilisant les d\'efinitions de $\cA$ et $\cB$, on trouve
\begin{eqnarray*}
\log \cG_{\beta,\alpha}
& \geq & -(D-1) \min \left( \log \left( \length(\cG_{1}) \right), \log \left( \length(\cG_{2}) \right) \right) \\
&      & +(K-1) \log (\max(1, |\beta|)) + (L-1) \log (\max(1, |\alpha|)) \\
&      & - D(K-1)\hgt(\beta) - D(L-1)\hgt(\alpha) \\
& \geq & -(D-1) \log \left( \min \left( d_{\mu}^{K-1}, \mu! \right) H \left( \cM_{0} \right)
                2^{\mu+K-1} e^{L-1}\sqrt{L} \right) \\
&      & -(K-1) \log (\cB) - (L-1)\log (\cA).
\end{eqnarray*}
\hfill $\Box$
\end{dem}

\subsection{Majoration de $\cG_{\beta,\alpha}$}

Pour v\'erifier les hypoth\`eses du lemme~\ref{lem:4.10}, il faut conna\^itre
une valeur admissible du param\`etre $\cN$. Dans ce but, nous allons
\'etablir la majoration analytique que voici.

\begin{lem}
\label{lem:4.16}
Soit $E$ un nombre r\'eel sup\'erieur \`a $1$. Soient $K$ et $L$ des nombres
entiers strictement positifs avec $L \geq 2$. On pose
\begin{eqnarray*}
\cN
& = & \max \left\{ E|\beta|L, \log(L-1) + (L-1)\log \left( e^{|\beta|}+\epsilon \right) \right\} \\
&   & + \log ((K-1)!) + (\mu+1) \log(L+1) + \log(2).
\end{eqnarray*}
Alors
$$
\max \left\{ \sum_{k=0}^{K-1} \sum_{\ell=0}^{L-1} |\Phi_{k,\ell}^{(\mu)}(z\beta)|_{E},
\sum_{k=0}^{K-1} \sum_{\ell=0}^{L-1} |w_{k,\ell}| \right\}
\leq e^{\cN}.
$$
\end{lem}

\begin{dem}
On estime $|\Phi_{k,\ell}^{(\mu)}(z\beta)|_{E}$.
On a 
$$
\left| \sum_{j=0}^{\min(\mu,k)} \frac{k!z^{k-j}}{(k-j)!} \right|
\leq k! e^{|z|}.
$$

Ainsi,
\begin{eqnarray*}
\left|\sum_{j=0}^{\min(\mu,k)} {\mu \choose j} k!\ell^{\mu-j} \frac{z^{k-j}}{(k-j)!} \right|
& \leq & \left| \sum_{j=0}^{\min(\mu,k)} {\mu \choose j} \ell^{\mu-j} \right|
         \cdot \left| \sum_{j=0}^{\min(\mu,k)} \frac{k!z^{k-j}}{(k-j)!} \right| \\
& \leq & (\ell+1)^{\mu} k! e^{|z|}
\end{eqnarray*}
et, d'apr\`es la d\'efinition de $\Phi_{k,\ell}^{(\mu)}(z)$ en Section 1.2,
\begin{eqnarray*}
\left| \Phi_{k,\ell}^{(\mu)}(z) \right|
&   =  & \left| \delta^{\mu} \left( X^{k}Y^{\ell} \right) (z,e^{z}) \right|
    =    \left| e^{\ell z} \sum_{j=0}^{\min(\mu,k)} {\mu \choose j}
         \frac{k!}{(k-j)!}z^{k-j}\ell^{\mu-j} \right| \\
& \leq & e^{\ell|z|} (\ell+1)^{\mu} k! e^{|z|}.
\end{eqnarray*}

Par cons\'equent
$$
|\Phi_{k,\ell}^{(\mu)}(z\beta)|_{E}
\leq e^{\ell|\beta z|} (\ell+1)^{\mu} k! e^{|\beta z|}
\leq e^{(\ell+1)|\beta z|} (\ell+1)^{\mu} k!.
$$

Et on en d\'eduit 
\begin{eqnarray*}
\sum_{k=0}^{K-1} \sum_{\ell=0}^{L-1} |\Phi_{k,\ell}^{(\mu)}(z\beta)|_{E}
& \leq &
e^{L|\beta E|} \sum_{k=0}^{K-1} k! \sum_{\ell=1}^{L} \ell^{\mu} \\
& < &
e^{L|\beta E|}\frac{(L+1)^{\mu+1}}{\mu+1} \sum_{k=0}^{K-1} k! \\
& < &
e^{L|\beta E|}\frac{(L+1)^{\mu+1}}{\mu+1} 2(K-1)!.
\end{eqnarray*}

Pour $\ell \geq 0$, on a 
\begin{eqnarray*}
\left| w_{k,\ell} \right|
& = & \left| \frac{\alpha^{\ell}-e^{\beta\ell}}{e^{\beta}-\alpha}
\sum_{j=0}^{\min(\mu,k)} {\mu \choose j} \frac{k!\ell^{\mu-j}}{(k-j)!} \beta^{k-j} \right| \\
& \leq & \left| \frac{\alpha^{\ell}-e^{\beta\ell}}{e^{\beta}-\alpha} \right| e^{|\beta|} k! (\ell+1)^{\mu} \\
& \leq & \ell\max(\alpha,e^{|\beta|})^{\ell-1} e^{|\beta|} k! (\ell+1)^{\mu} \\
& \leq & \ell(e^{|\beta|}+\epsilon)^{\ell} k! (\ell+1)^{\mu}.
\end{eqnarray*}

Ainsi, comme ci-dessus,
\begin{eqnarray*}
\sum_{k=0}^{K-1} \sum_{\ell=0}^{L-1} \left| w_{k,\ell} \right|
& \leq &
(L-1) (e^{|\beta|}+\epsilon)^{L-1} \sum_{k=0}^{K-1} k! \sum_{\ell=1}^{L} \ell^{\mu} \\
& < &
(L-1) (e^{|\beta|}+\epsilon)^{L-1} \frac{(L+1)^{\mu+1}}{\mu+1} \sum_{k=0}^{K-1} k! \\
& < &
(L-1) (e^{|\beta|}+\epsilon)^{L-1}\frac{(L+1)^{\mu+1}}{\mu+1} 2(K-1)!.
\end{eqnarray*}

%

D'o\`u le lemme~\ref{lem:4.16}. \hfill$\Box$
\end{dem}

\vspace{3.0mm}

Ainsi, sous les hypoth\`eses du lemme~\ref{lem:4.10}, on obtient la proposition
suivante.

\begin{prop}
\label{prop:4.17alg}
Soient $E>1$ un nombre r\'eel, $K$ et $L$ deux entiers strictement
positifs avec $L \geq 2$. On suppose $\epsilon<E^{-KL}$.

On a 
\begin{eqnarray*}
\log \cG_{\beta,\alpha}
& \leq & -(KL-\mu-1)\log (E) + \log(4) + \log ((K-1)!) \\
&      & + \max \left\{ E|\beta|L, \log(L)+L\log\left(e^{|\beta|}+\epsilon \right) \right\}
+ (\mu+1)\log(L+1) \\
&      & + \log \left( \min \left( d_{\mu}^{K-1}, \mu! \right) \right)
         + \log \left( H(\cM_{0}) \right) - \log (\mu!).
\end{eqnarray*}
\end{prop}

\begin{dem}
La proposition~\ref{prop:4.17alg} est une cons\'equence des lemmes~\ref{lem:4.10}
et~\ref{lem:4.16} et de la relation entre $\cD$ et $\cG_{\beta,\alpha}$ donn\'ee
par (\ref{eq:defG}).
\hfill$\Box$
\end{dem}

\subsection{Fin de la d\'emonstration du th\'eor\`eme~\ref{theo:4.6bisalg}}

On raisonne par l'absurde en supposant $\epsilon \leq E^{-KL}$. On va utiliser
la minoration de $\cG_{\beta,\alpha}$ fournie par la proposition~\ref{prop:trivalg}
et la majoration donn\'ee par la proposition~\ref{prop:4.17alg}.

Comme $L \geq 2$, on majore $\mu$ par $L-2$ et $\mu\log(L)-\log(\mu!)$ par
$L-0.5\log(2\pi L)$, en utilisant le lemme~\ref{lem:factorials}(i) et en faisant
quelques calculs pour de petites valeurs de $L$. Donc on obtient
\begin{eqnarray*}
&      & -D \log \left( H \left( \cM_{0} \right) \min \left( d_{L-2}^{K-1}, (L-2)! \right) \right)
         -(D-1) \log \left( 2^{L+K-3} e^{L-1}\sqrt{L} \right) \\
&      & -(K-1) \log (\cB) - (L-1)\log (\cA) \\
& \leq & -(KL-L+1)\log (E) + \log(4) + \log ((K-1)!) + L - \frac{\log(2\pi L)}{2} \\
&      & + \max \left\{ E|\beta|L, \log(L-1)+(L-1)\log\left(e^{|\beta|}+E^{-KL} \right) \right\}.
\end{eqnarray*}

Par ailleurs, on peut supposer $KL\log (E) \geq DLK\log (2) + L \log(E)$,
soit $(K-1)L \log (E) \geq DLK\log(2)$. Puisque $K$ et $L$ sont des entiers
strictement positifs, le membre de droite est positif et donc $K \geq 2$.
On a alors
\begin{eqnarray*}
(L-1)\log\left(e^{|\beta|}+E^{-KL}\right)
& \leq & |\beta|(L-1)+(L-1)\log \left(1+2^{-2L} \right) \\
&  <   & |\beta|L+L2^{-L} \leq |\beta|L+0.061.
\end{eqnarray*}

Donc
$$
\max \left\{ E|\beta|L, \log(L-1)+(L-1)\log\left(e^{|\beta|}+E^{-KL}\right) \right\}
< E|\beta|L + \log(L-1) + 0.061.
$$

En appliquant ceci avec la proposition~\ref{prop:3.10} et $\log(4)-\log(2\pi)/2+0.061<0.53$,
on trouve
\begin{eqnarray*}
&   & -D \log \left\{ D_{K-1,L-1} \left( \frac{\sqrt{3}}{2} e \frac{\min(K,L)}{\sqrt{L}}d_{L-1} \right)^{K-1}
      2^{KL+L} \frac{\sqrt{6}}{16L} \min \left( d_{L-2}^{K-1}, (L-2)! \right) \right\} \\
&   & -(D-1) \log \left( 2^{L+K-3} e^{L-1}\sqrt{L} \right)
      -(K-1) \log (\cB) - (L-1)\log (\cA) \\
& < & -(KL-L+1)\log (E) + \log ((K-1)!) + L + E|\beta|L + 0.5\log(L) + 0.53.
\end{eqnarray*}

En regroupant et r\'earrangeant les diff\'erents termes, on obtient
\begin{eqnarray}
&   & KL \log (E) \\
& < & DKL \log (2) + D(K-1) \log \left( \sqrt{3} e \frac{\min(K,L)}{\sqrt{L}}d_{L-1} \right) \nonumber \\
&   & +D \log \left\{ (4e)^{L-1} \min \left( d_{L-2}^{K-1}, (L-2)! \right) \right\}
      +D \log \left( \frac{D_{K-1,L-1}}{\sqrt{L}} \right) \nonumber \\
&   & + \log ((K-1)!) +(K-1) \log (\cB/2) + (L-1)\log (\cA/2) \nonumber \\
&   & + E|\beta|L +(L-1)\log (E) + \log \left( 2 e \right) + D \log \left( \frac{\sqrt{6}}{16} \right) + 0.53. \nonumber
\end{eqnarray}

La majoration $D\log \left( \sqrt{6}/(16\sqrt{L}) \right) + \log(2e)+0.53 < 0$
permet de conclure cette d\'emonstration du th\'eor\`eme~\ref{theo:4.6bisalg}.

\section{D\'emonstration du corollaire~\ref{coro:4}}

On choisit $K=\lfloor c_{1}|\beta| \rfloor$ et $L=\lfloor c_{2} \log(|\beta|) \rfloor$
avec des constantes r\'eelles strictement positives $c_{1}$ et $c_{2}$. Nous
d\'eterminerons les valeurs optimales pour pour $c_{1}$, $c_{2}$ et $E$. Cela
nous permettra de d\' emontrer le corollaire~\ref{coro:4}.

On peut \'{e}crire le membre de gauche de l'in\'egalit\'e (\ref{eq:theoTalg})
de la fa\c{c}on suivante
$$
c_{1}c_{2} |\beta|\log(|\beta|) \log (E) + o(|\beta|\log(|\beta|)).
$$

Puisque les nombres $\alpha$ et $\beta$ sont des entiers alg\'{e}briques dans
un corps imaginaire quadratique, on a $D=1$ et on peut prendre $\cA=\cB=1$.
On majore $D_{K-1,L-1}$ par $(K-1)!$, comme le logarithme du d\'enominateur
de $D_{K-1,L-1}$ est en $o(|\beta| \log (|\beta|))$. Aussi on \'ecrit $d_{x}$
sous la forme $\exp \left( \psi(x) \right)$, o\`u $\psi(x)=\sum_{p^{k}\leq x} \log p$
et on remplace $x!$ par $\Gamma(x+1)$. \'Etant donn\'e que $e\sqrt{3L}$ et
$(4e)^{L-1}$ sont tous deux en $o(|\beta| \log (|\beta|))$, on peut \'ecrire le
membre de droite de l'in\'egalit\'e (\ref{eq:theoTalg}) sous la forme
\begin{eqnarray*}
& & c_{1}c_{2}|\beta| \log(|\beta|)\log 2 + c_{1}|\beta|\psi \left( c_{2} \log(|\beta|) \right)
    + 2\log ( \Gamma(c_{1}|\beta|)) \\
& & + \min \left( c_{1}|\beta| \psi \left( c_{2} \log(|\beta|) \right), \log \left( \Gamma(c_{2} \log(|\beta|)) \right) \right) \\
& & + E|\beta| \left( c_{2} \log(|\beta|) \right)
    + c_{2} \log(|\beta|)\log E + o(|\beta| \log (|\beta|)).
\end{eqnarray*}

Des d\'eveloppements asymptotiques $\log \Gamma(x)=x\log(x)+o(x\log x)$ et $\psi(x)=x+o(x)$,
on d\' eduit que le membre de droite de l'in\'egalit\'e (\ref{eq:theoTalg}) est
$$
c_{1}c_{2} (1+\log (2)) |\beta| \log (|\beta|) + 2c_{1}|\beta| \log(|\beta|)
+ Ec_{2} |\beta|\log(|\beta|) + o(|\beta|\log(|\beta|)).
$$

En combinant ces expressions pour les membres de gauche et droite de l'in\'egalit\'e
(\ref{eq:theoTalg}) et en divisant ces deux m\^emes membres par $|\beta|\log(|\beta|)$,
la condition \`a satisfaire devient 
\begin{equation}
\label{eq:coro4-1}
c_{1}c_{2} \log E > c_{1}c_{2} \left( 1 + \log (2) \right) + 2c_{1} + Ec_{2} + o(1).
\end{equation}
En ignorant le terme $o(1)$, pour le moment, et en r\'esolvant l'\'equation 
$$
c_{1}c_{2} \log E = c_{1}c_{2} \left( 1 + \log (2) \right) + 2c_{1} + Ec_{2}
$$
pour $c_{1}$, on obtient
$$
c_{1} = \frac{c_{2}E}{c_{2}(\log(E)-1-\log(2))-2},
$$
qui s'\'ecrit 
\begin{equation}
\label{eq:coro4-2}
c_{1}c_{2} \log E
= \frac{c_{2}^{2}E\log(E)}{c_{2}(\log(E)-1-\log(2))-2}.
\end{equation}

En d\'erivant cette fonction de $c_{2}$ dans le membre de droite et en r\'esolvant
l'\'equation permettant de trouver les points annulant cette d\'eriv\'ee, on
trouve la valeur
$$
c_{2} = \frac{4}{\log(E)-1-\log(2)}.
$$

Pour cette valeur de $c_{2}$, le membre de droite de l'\'egalit\'e~(\ref{eq:coro4-2})
vaut
$$
\frac{8E\log(E)}{(\log(E)-1-\log(2))^{2}}.
$$

En d\'erivant cette expression par rapport \`a $E$ et en r\'esolvant \`a nouveau
l'\'equation permettant de trouver les points annulant cette nouvelle d\'eriv\'ee,
on trouve la valeur
$$
E=\exp \left( 1+\frac{\log(2)}{2}+\frac{1}{2}\sqrt{\log^{2}(2)+8\log(2)+8} \right)
=25.0059\ldots
$$
En posant $\gamma=\log^{2}(2)+8\log(2)+8$, on obtient
\begin{eqnarray*}
c_{1}c_{2}\log(E)
& = & \frac{16\sqrt{2}\exp \left( 1+(1/2)\sqrt{\gamma} \right)
      \left( 2+\log(2) + \sqrt{\gamma} \right)}
      {\left( \sqrt{\gamma} - \log(2) \right)^{2}}
=276.55\ldots
\end{eqnarray*}
Il r\'esulte de cette \'etude que pour tout $\epsilon>0$, il existe $B_{0}$ tel
que pour $|\beta| > B_{0}$, l'in\'egalit\'e (\ref{eq:coro4-1}) soit v\'erifi\'ee
avec $c_{1}c_{2}\log(E)=276.55\ldots+\epsilon$. D'o\`u le corollaire~\ref{coro:4}.

\noindent
{\sc S. Kh\'{e}mira\\
Paris, France\\
Courriel: khemira@math.jussieu.fr}

\vspace{3.0mm}

\noindent
{\sc P. Voutier\\
London, UK\\
Courriel: paul.voutier@gmail.com}
\end{document}